%% file: m3-I-1.tex
\input m3-macs

\pageno=5

\tinfo I.1.5-18

\SetTFLinebox{\gtp }
\SetFLinebox{\gtv3 }
\SetHLinebox{\issn}

\H 1.  Higher dimensional local fields

Igor Zhukov

\SetAuthorHead{I. Zhukov}
\SetTitleHead{Part I. Section 1. Higher dimensional local fields
 \qquad\qquad}

We give here basic definitions related to $n$-dimensional local fields. For detailed
exposition, see \cite{P} in the equal characteristic case,
\cite{K1, \S8} for the two-dimensional case and \cite{MZ1}, \cite{MZ2} 
for the general case.
Several properties of the topology on the multiplicative group are discussed
in \cite{F}.

\HH 1.1. Main definitions

Suppose that we are given a surface $S$ over a finite field
of characteristic $p$, a curve $C\subset S$, and a point $x\in C$
such that both $S$ and $C$ are regular at $x$. Then
one can attach to these data  the quotient field
of the completion $\widehat{(\widehat{\Cal O}_{S,x})}_C$ of the localization at $C$   of the completion $\widehat{\Cal O}_{S,x}$
of the local ring $\Cal O_{S,x}$ of $S$ at $x$.
This is a two-dimensional local field over a finite field, i.e., a \cdvf
with local residue field. More generally, an $n$-dimensional  local field $F$ is
a \cdvf with $(n-1)$-dimensional residue field.
(Finite fields are considered as 0-dimensional local fields.)

 \df  Definition

A \cdvf $K$ is said to have the structure of an $n$-dimensional local
field if there is a chain of fields $K=K_{n},K_{n-1},\dots,K_{1},K_{0}$
where $K_{i+1}$ is a  \cdvf  with residue field $K_{i}$ and $K_{0}$ is a finite
field.
The field $k_K=K_{n-1}$ (resp.\ $K_0$) is said to be the {\it first} (resp.\ the {\it last})
residue field of $K$.
\enddf

\rk Remark

Most of the properties of $n$-dimensional local fields do not change
if one requires that the last residue $K_0$ is {\it perfect} rather than
finite. To specify the exact meaning of the word, $K$ can be referred
to as an $n$-dimensional local field {\it over} a finite (resp.\ perfect)
field. One can consider an $n$-dimensional local field over an arbitrary
field $K_0$ as well. However, in this volume mostly the higher
local fields over finite fields are considered.
\endrk

\eg Examples

1. $ \Bbb F_{q}(( X_{1}) ) \dots((
X_{n})) $.
2.  $k( ( X_{1})) \dots((
X_{n-1}))$, $k$ a finite extension of $\Bbb Q_p$.

3. For a  \cdvf  $F$ let
$$
K=F\left\{ \!\left\{ T\right\}\! \right\} =\biggl\{ \sum_{-\infty}^{+\infty
}a_{i}T^{i} : a_{i}\in F,\;\innf v_F(a_{i})>-\infty,\;\lim\limits_{i\rightarrow-\infty
}v_F(a_{i})=+\infty \biggr\} .
$$
Define
$v_{K}(\sum a_{i}T^{i})=\minn v_F(a_{i})$.
Then $K$ is a  \cdvf  with residue field $k_{F}\left( \left( t\right) \right) $.

Hence for a local field $k$
the fields
$$k\left\{\! \left\{ T_{1}\right\}\! \right\} \dots\left\{\! \left\{ T_{m}\right\}\!
\right\} ( ( T_{m+2}) ) \dots((T_{n})), \quad 0\le m\le n-1 $$
are  $n$-dimensional local fields
(they are called {\it standard fields}).
\endeg

\rk Remark

 $K\left( \left( X\right) \right) \left\{\! \left\{ Y\right\}\!
\right\} $ is isomorphic to $K\left( \left( Y\right) \right) \left( \left(
X\right) \right) $.
\endrk

 \df  Definition

An $n$-tuple of elements $t_{1},\dots,t_{n}\in K$ is
called a {\it system of local parameters} of $K$,
if $t_{n}$ is a prime element of $K_{n}$, $%
t_{n-1}$ is a unit in $ \Cal O _{K}$ but its residue in $K_{n-1}$ is a
prime element  of $K_{n-1}$, and so on.
\enddf

For example, for $K=k\left\{\! \left\{ T_{1}\right\}\! \right\} \dots\left\{\! \left\{ T_{m}\right\}\!
\right\} ( ( T_{m+2}) ) \dots((T_{n}))$, a convenient system of local parameter is
$T_1,\dots,T_m,\pi,T_{m+2},\dots,T_n$, where $\pi$ is a prime element of $k$.

Consider the maximal $m$ such that $\chr(K_{m})=p$;
we have $0\le m\le n$. Thus,  there are
$n+1$ types of $n$-dimensional local fields:
fields of characteristic $p$ and
fields with $\chr(K_{m+1})=0$, $\chr(K_{m})=p$, $0\le m\le n-1$.
Thus, the mixed characteristic case is the case $m=n-1$.

Suppose that $\chr( k_K)=p$, i.e., the above $m$ equals either $n-1$ or $n$. Then
the set of Teichm\"uller representatives $\Cal R$ in $\Cal O_K$ is a field
isomorphic to $K_0$.

\th Classification Theorem

 Let $K$ be an $n$-dimensional local field.
Then
\Roster

\Item{(1)} $K$ is isomorphic to
$ \Bbb F_{q}( ( X_{1}) )
\dots( ( X_{n}) ) $
if $\chr(K)=p$;

\Item{(2)} $K$ is isomorphic to
$k((X_{1}))\dots(( X_{n-1}))$, $k$ is a local field,
 if  $\chr(K_{1})=0$;

\Item{(3)} $K$ is a finite extension of a standard field
$
k\left\{\! \left\{ T_{1}\right\}\! \right\} \dots\left\{\! \left\{ T_{m}\right\}
\!\right\} ( ( T_{m+2}) ) \dots((
T_{n}))$ and there is a finite extension of $K$ which is
a standard field
if $\chr(K_{m+1})=0$, $\chr(K_{m})=p$.
\endRoster
\endth

\pf Proof

In the equal characteristic case the statements follow from
the well known classification theorem
for complete discrete valuation fields of equal characteristic.
In the mixed characteristic case
let
$k_{0}$ be the fraction field of $W( \Bbb F_{q})$ and let
$T_{1},...,T_{n-1}, \pi$ be a
system of local parameters of $K$.
Put
$$K'=k_{0}\left\{\! \left\{ T_{1}\right\}\!
\right\} \dots \left\{\! \left\{ T_{n-1}\right\}\! \right\}.$$
Then $K'$ is an absolutely unramified complete discrete valuation field, and the (first) residue fields
of $K'$ and $K$ coincide.
Therefore, $K$ can be viewed as a finite extension of $K'$ by \cite{FV, II.5.6}.

Alternatively, let $t_{1},\dots ,t_{n-1}$ be any liftings of a system of local parameters of $k_K$.
Using the canonical lifting $h_{t_1,\dots,t_{n-1}}$ defined below, one can construct an embedding $K'\hookrightarrow K$
which identifies $T_i$ with $t_i$.

To prove the last assertion of the theorem, one can use {\it Epp's theorem on
elimination of wild ramification}
(see 17.1) which asserts that
there is a finite extension $l/k_0$ such
that $e\left( lK/lK'\right) =1$.  Then $lK'$ is standard and $lK$
is standard, so $K$ is a subfield of $lK$. See \cite{Z} or \cite{KZ} for details
and a stronger statement.
\qed\endpf

\df Definition

The {\it lexicographic  order} of $ \Bbb Z^{n}$:
$
{\bold i}=(i_{1},\dots,i_{n})\le {\bold j}=(j_{1},\dots,j_{n})$
if and only if
$$ i_{l}\le
j_{l},\;i_{l+1}=j_{l+1},\dots,i_{n}=j_{n} \text{ for some $l \le n$}.
$$

Introduce $\bold v=(v_{1},\dots,v_{n})\colon K^{*}\rightarrow \Bbb Z^{n}$
as
$v_{n}=v_{K_{n}}$, $v_{n-1}(\alpha)  =v_{K_{n-1}}(\alpha_{n-1})$
where $\alpha_{n-1}$ is the residue of $\alpha t_{n}^{-v_{n}(\alpha)}$ in $K_{n-1}$%
, and so on. The map $\bold v$ is a valuation;
this is a so called discrete valuation of rank $n$.
Observe that for $n>1$ the valuation $\bold v$ does depend on the choice of $%
t_{2},\dots,t_{n}$. However, all the valuations obtained this way
are in the same class of equivalent valuations.
\enddf

Now we define several objects which do not depend on the choice
of a system of local parameters.

 \df  Definition

{}

$ {\CO} _{K}=\left\{ \alpha\in K: {\bold v}  %
(\alpha)\ge 0\right\} $, $ {\CM}_{K}=\left\{ \alpha\in K: {\bold v}  %
(\alpha)>0\right\} $, so $ {\CO} _{K}/ {\CM}_{K}\simeq K_{0}$.
The group of {\it principal units} of $K$ with respect
to the valuation $\bold v$ is
$V_K=1+{\CM}_K$.

\enddf

 \df  Definition

{}
$$ P\left( i_{l},\dots,i_{n}\right) =P_K\left( i_{l},\dots,i_{n}\right)=\{
\alpha\in K: \left( v_{l}(\alpha),\dots,v_{n}(\alpha)\right) \ge\left(
i_{l},\dots ,i_{n}\right) \} .$$
\enddf

In particular, $\CO_K=P\bigl(\underset{n}\to{\underbrace{0,\dots,0}}\,\bigr)$,
$\CM_K=P\bigl(1,\underset{n-1}\to{\underbrace{0,\dots,0}}\,\bigr)$,
whereas $\Cal O_K=P(0)$, $\Cal M_K=P(1)$.
Note that if $n>1$, then 
 $$\cap_i M_K^i=P\bigl(1,\underset{n-2}\to{\underbrace{0,\dots,0}}\,\bigr),$$ 
since 
 $t_2=t_1^{i-1}(t_2/t_1^{i-1})$.

\th Lemma

 The set of all non-zero ideals of $ {\CO} _{K}$ consists of all
$$\{P \left(
i_{l},\dots ,i_{n}\right) : \left( i_{l},\dots
,i_{n}\right)  \ge \left( 0,\dots ,0\right), \quad
1 \le l \le n\}.$$
The ring  ${\CO}_{K}$ is not Noetherian for $n>1$.
\endth

\pf Proof

Let $J$ be a non-zero ideal of $\CO_K$. Put $i_n=\min\{v_n(\alpha):\alpha\in J\}$.
If $J=P\left(i_n\right)$, then we are done. Otherwise, it is clear that
$$i_{n-1}:=\inf\{v_{n-1}(\alpha): \alpha\in J,v_n(\alpha)=i_n\}>-\infty.$$
If $i_n=0$, then obviously $i_{n-1}\ge0$.
Continuing this way, we construct $(i_l,\dots,i_n)\ge(0,\dots,0)$, where either $l=1$
or
$$
i_{l-1}=\inf\{v_{l-1}(\alpha) : \alpha\in J,v_n(\alpha)=i_n,\dots,v_l(\alpha)=i_l\}=-\infty.
$$
In both cases it is clear that $J=P \left(i_{l},\dots ,i_{n}\right)$.

The second  statement is immediate from $P(0,1)\subset P(-1,1)\subset P(-2,1) \dots$.
\qed\endpf

For more on ideals in $O_K$ see subsection~3.0 of Part~II.

\HH 1.2. Extensions

Let $L/K$ be a finite extension.  If $K$ is an
$n$-dimensional local field, then so is $L$.

 \df  Definition

Let
 $ t_{1},\dots,t_{n}$  be a system of local parameters of $K$
and let
$t_{1}^{^{\prime}},\dots,t_{n}^{^{\prime}}$ be a system of local parameters of $L$.
Let
$ \bold{v}  , \bold{v}  ^{\prime}$ be the corresponding valuations.
Put
$$
E(L|K):=\bigl( v_{j}^{^{\prime}}(t_{i})\bigr) _{i,j}=\left(
\matrix e_{1} & 0 &\dots &0 \\
\dots & e_2&\dots & 0 \\
\dots&\dots&\dots & 0 \\
\dots & \dots &\dots & e_{n}
\endmatrix
\right) ,
$$
where $e_i=e_{i}(L|K)=e(L_{i}|K_{i})$, $i=1,\dots,n$.
Then  $e_{i}$ do not depend on the choice of parameters,
and
$\left| L:K\right| =f(L|K)\prod_{i=1}^{n} e_{i}(L|K) $, where $f(L|K)=\left| L_{0}:K_{0}\right|$ .
\enddf

The expression ``unramified extension'' can be used for extensions $L/K$ with

\noindent $e_n(L|K)=1$
and $L_{n-1}/K_{n-1}$ separable. It can be also used in a narrower sense, namely,
for extensions $L/K$ with $\prod_{i=1}^{n} e_{i}(L|K)=1$. To avoid ambiguity, sometimes
one speaks of a ``semiramified extension'' in the former case and a ``purely unramified
extension'' in the latter case.

\HH 1.3. Topology on $K$


Consider an example of $n$-dimensional local field
$$K=k\left\{\! \left\{ T_{1}\right\}\! \right\} \dots\left\{\! \left\{ T_{m}\right\}\!
\right\} ( ( T_{m+2}) ) \dots((T_{n})).$$
Expanding elements of $k$ into power series in $\pi$ with coefficients in $\Cal R_k$, one can write elements of $K$ as formal power series in $n$ parameters. To make them convergent power series we should introduce a topology in $K$
which takes into account topologies of the residue fields. We do not make $K$ a topological field this way,
since  multiplication is only sequentially continuous in this topology.
However, for class field theory sequential continuity seems to be more important than continuity.

\HHH 1.3.1

\df Definition

\Roster

\Item{(a)}  If $F$
has a topology, consider the following topology on $K=F\left( \left( X\right) \right) $.
For a sequence of neighbourhoods of zero $\left( U_{i}\right) _{i\in \Bbb Z}$ in $F$, $%
U_{i}=F$ for $i\gg 0$, denote $U_{\left\{ U_{i}\right\} }=\left\{ \sum
a_{i}X^{i}:a_{i}\in U_{i}\right\} $. Then all $U_{\left\{ U_{i}\right\} }$ constitute a base of open neighbourhoods of 0 in
$F\left( \left( X\right) \right) $.
In particular, a sequence $u^{(n)}=\sum a_i^{(n)} X^i$ tends to 0
if and only if there is
an integer $m$ such that $u^{(n)}\in X^mF[[X]]$ for all $n$
and the sequences $a_i^{(n)}$ tend to 0 for every $i$.
\smallskip
\Item{{}}
Starting with the discrete topology on the last residue field, this construction is used
to obtain a well-defined topology on an $n$-dimensional local field of characteristic $p$.

\Item{(b)}  Let $K_n$ be of mixed characteristic.  Choose a system of local parameters $t_1,\dots,t_n$ $=\pi$ of $K$.
The choice of $t_1,\dots,t_{n-1}$ determines a canonical lifting $$h=h_{t_1,\dots,t_{n-1}}\:K_{n-1}\to \Cal O_K$$
(see below).
Let $\left( U_{i}\right) _{i\in \Bbb Z}$ be a system of neighbourhoods of zero in $K_{n-1}$, $U_{i}=K_{n-1}$ for $%
i\gg 0 $.
Take the system of all $U_{\left\{ U_{i}\right\} }=\left\{ \sum
h(a_{i})\pi^{i},\;a_{i}\in U_{i}\right\} $ as a base of open neighbourhoods
of 0 in $K$.
This topology is well defined.

\Item{(c)} In the case $\chr (K)=\chr (K_{n-1})=0$ we apply constructions (a) and (b)  to obtain a topology on $K$
which depends  on  the choice of the coefficient
subfield of $K_{n-1}$ in $\Cal O_K$.

\endRoster
\enddf

The definition of the canonical lifting $h_{t_1,\dots,t_{n-1}}$ is rather
complicated. In fact, it is worthwhile to define it for any $(n-1)$-tuple
$(t_1,\dots,t_{n-1})$ such that $v_i(t_i)>0$ and $v_j(t_i)=0$ for $i<j\le n$.
We shall give an outline of this construction, and the details can be found
in \cite{MZ1, \S1}.

Let $F=K_0((\Overline{t_1}))\dots((\Overline{t_{n-1}}))\subset K_{n-1}$.
By a lifting we mean a map $h\:F\to \Cal O_K$ such that the residue of $h(a)$
coincides with $a$ for any $a\in F$.

Step 1. An auxiliary lifting $H_{t_1,\dots,t_{n-1}}$ is uniquely determined
by the condition
$$
\split
&H_{t_1,\dots,t_{n-1}} \biggl(\sum_{i_1=0}^{p-1}\dots\sum_{i_{n-1}=0}^{p-1}
  {\Overline{t_1}}^{i_1}\dots  {\Overline{t_{n-1}}}^{i_{n-1}}a_{i_1,\dots,i_{n-1}}^p\biggr)
\\
&\qquad =\sum_{i_1=0}^{p-1}\dots\sum_{i_{n-1}=0}^{p-1}
  t_1^{i_1}\dots t_{n-1}^{i_{n-1}} (H_{t_1,\dots,t_{n-1}}(a_{i_1,\dots,i_{n-1}}))^p.
\endsplit
$$

Step 2. Let $k_0$ be the fraction field of $W(K_0)$. Then
$K'=k_0\{\!\{T_1\}\!\}\dots\allowmathbreak\{\!\{T_{n-1}\}\!\}$
is an $n$-dimensional local field with the residue field $F$.
Comparing the lifting $H=H_{T_1,\dots,T_{n-1}}$ with the lifting $h$
defined by
$$
h\bigl(\sum_{{{{\bold r}}}\in\Bbb Z^{n-1}}\theta_{{{\bold r}}}{ {\Overline{T_1}}^{r_1}\dots  {\Overline{T_{n-1}}}^{r_{n-1}}}\bigr)
=
\sum_{{{{\bold r}}}\in\Bbb Z^{n-1}} [\theta_{{{\bold r}}}]{T_1^{r_{1}}\dots T_{n-1}^{r_{n-1}}},
$$
we introduce the maps $\lambda_i\:F@>>>F$ by the formula
$$
h(a)=H(a)+pH(\lambda_1(a))+p^2H(\lambda_2(a))+\dots
$$

Step 3. Introduce $h_{t_1,\dots,t_{n-1}}\:F@>>>\Cal O_K$ by the formula
$$
h_{t_1,\dots,t_{n-1}}(a)=H_{t_1,\dots,t_{n-1}}(a)+pH_{t_1,\dots,t_{n-1}}(\lambda_1(a))+p^2H_{t_1,\dots,t_{n-1}}(\lambda_2(a))+\dots.
$$

\rk Remarks

1. Observe that for a standard field
 $K=k\left\{\!\left\{ T_{1}\right\}\! \right\}\!
\dots\left\{ \!\left\{ T_{n-1}\right\}\! \right\}, $ we have 
$$
h_{T_1,\dots,T_{n-1}}\:  \sum \theta_{{\bold i}}\Overline{T_{1}}^{i_{1}}\dots  \Overline{T_{n-1}}^{i_{n-1}}
\mapsto\sum [\theta_{{\bold i}}]T_{1}^{i_{1}}\dots T_{n-1}^{i_{n-1}},
$$
where $\Overline{T_j}$ is the residue of $T_j$ in $k_K$, $j=1,\dots,n-1$.

2. The idea of the above construction is to find
a field $k_0\{\!\{t_1\}\!\}\dots\{\!\{t_{n-1}\}\!\}$
isomorphic to $K'$ inside $K$
without a priori given topologies on $K$ and $K'$. More precisely, let
$t_1,\dots,t_{n-1}$ be as above.  For $a=\sum_{-\infty}^\infty p^i h(a_i) \in K'$, let
$$
f_{t_1,\dots,t_{n-1}}(a)=\sum_{-\infty}^\infty p^i h_{t_1,\dots,t_{n-1}}(a_i)
$$
Then $f_{t_1,\dots,t_{n-1}}\:K'@>>>K$ is an embedding of $n$-dimensional
complete fields such that
$$
f_{t_1,\dots,t_{n-1}}(T_j)=t_j, \qquad j=1,\dots,n-1
$$
(see \cite{MZ1, Prop. 1.1}).

3. In the case of a standard mixed characteristic field the following alternative
construction of the same topology is very useful.

Let $K=E\{\!\{X\}\!\}$, where $E$ is an $(n-1)$-dimensional local field; assume
that the topology of $E$ is already defined.
Let $\{V_i\}_{i\in\Bbb Z}$ be a sequence of neighbourhoods of zero
in $E$ such that

\noindent(i) there is $c\in\Bbb Z$ such that $P_E(c)\subset V_i$ for all $i\in \Bbb Z$;

\noindent(ii)  for every $l\in \Bbb Z$
we have  $P_E(l)\subset V_i$ for all sufficiently large $i$.

Put
$$
\Cal V_{\{V_i\}}=
\bigl\{\sum b_i
X^i :  b_i\in V_i \bigr\}.
$$
Then all the sets $\Cal V_{\{V_i\}}$ form a base of neighbourhoods of 0 in $K$.
(This is an easy but useful exercise in the 2-dimensional case;
in general, see Lemma 1.6 in \cite{MZ1}).

4. The formal construction of $h_{t_1,\dots,t_{n-1}}$ works also in  case $\chr (K)=p$, and one need not consider
this case separately. However, if one is interested in  equal characteristic case only,
all the treatment can be considerably simplified. (In fact, in this case $h_{t_1,\dots,t_{n-1}}$
is just the obvious embedding of $F\subset k_K$ into $\Cal O_K=k_K[[t_n]]$.)
\endrk

\HHH 1.3.2. Properties

\Roster

\Item{(1)} $K$ is a topological group which is complete and separated.

\Item{(2)} If $n>1$, then every base of neighbourhoods of 0
is uncountable. In
particular, there are maps which are sequentially continuous but not
continuous.

\Item{(3)} If $n>1$,  multiplication in $K$ is not continuous.  In fact,
$UU=K$ for every
open  subgroup $U$, since $U\supset P(c)$ for some $c$ and $U\not\subset P(s)$
for any $s$.
 However,  multiplication is sequentially
continuous:
$$
\alpha_{i}\rightarrow \alpha,\quad 0\not=\beta_{i}\rightarrow \beta \neq0 \Longrightarrow
\alpha_{i}\beta_{i}^{-1}\rightarrow \alpha\beta^{-1}.
$$

\Item{(4)}  The map
$K  \rightarrow K,\quad \alpha \mapsto c\alpha$ for $c\neq0$
is a homeomorphism.

\Item{(5)} For a finite extension $L/K$ the topology of $L=$ the topology of
finite dimensional vector spaces over $K$ (i.e., the product topology on $K^{|L:K|}$).
Using this property one can redefine the topology
first for ``standard'' fields
$$k\left\{\! \left\{ T_{1}\right\}\! \right\} \dots
\left\{\! \left\{ T_{m}\right\}\!\right\}\left( \left(
T_{m+2}\right) \right) \dots \left( \left( T_{n}\right) \right) $$ using the
canonical lifting $h$, and then for arbitrary fields as the topology of
finite dimensional vector spaces.

\Item{(6)} For a finite extension $L/K$ the topology of $K=$ the topology induced from $L$.
Therefore, one can use the Classification
Theorem and define the topology on $K$ as induced by that on $L$,
where $L$ is taken to be a standard $n$-dimensional local field.
\endRoster

\rk Remark

In practical work with higher local fields, both (5) and (6) enables one
to use the original definition of topology only in the simple case
of a standard field.

\HHH 1.3.3. About proofs

The outline of the proof of assertions in {1.3.1--1.3.2} is as follows.
(Here we concentrate on the most complicated case $\chr (K)=0$, $\chr (K_{n-1})=p$;
the case of $\chr(K)=p$ is similar and easier, for details see \cite{P}).

{\it Step 1} (see \cite{MZ1, \S1}). Fix first $n-1$ local parameters (or, more generally, any elements
$t_1,\dots,t_{n-1}\in K$ such that $v_i(t_i)>0$ and $v_j(t_i)=0$ for $j>i$).

Temporarily fix
$\pi_i\in K$ ($i\in\Bbb Z$), $v_n(\pi_i)=i$, and $e_j\in P_K(0)$,
$j=1,\dots,d$, so that $\{\Overline{e_j}\}_{j=1}^d$ is a basis
of the $F$-linear space $K_{n-1}$. (Here $F$ is as in {1.3.1}, and $\Overline\alpha$
denotes the residue of $\alpha$ in $K_{n-1}$.) Let $\{U_i\}_{i\in\Bbb Z}$ be a sequence
of neighbourhoods of zero in $F$, $U_i=F$ for all sufficiently large $i$. Put
$$
\Cal U_{\{U_i\}}=
\bigl\{\sum_{i\ge i_0}\pi_i \cdot \sum_{j=1}^d e_j h_{t_1,\dots,t_{n-1}}(a_{ij})
: a_{ij}\in U_i,i_0\in\Bbb Z\bigr\}.
$$
The collection of all such sets $\Cal U_{\{U_i\}}$ is denoted by $B_U$.

{\it Step 2} (\cite{MZ1, Th. 1.1}). In parallel one proves that

-- the set $B_U$ has a cofinal subset which consists of subgroups of $K$; thus, $B_U$ is
a  base of neighbourhoods of zero of a certain topological group $K_{t_1,\dots,t_{n-1}}$
with the underlying (additive) group $K$;

-- $K_{t_1,\dots,t_{n-1}}$ does not depend on the choice of $\{\pi_i\}$ and $\{e_j\}$;

-- property (4) in {1.3.2} is valid for $K_{t_1,\dots,t_{n-1}}$.

{\it Step 3} (\cite{MZ1, \S2}). Some properties of $K_{t_1,\dots,t_{n-1}}$ are established,
in particular, (1) in {1.3.2}, the sequential continuity of multiplication.

{\it Step 4} (\cite{MZ1, \S3}). The independence from
the choice of  $t_1,\dots,t_{n-1}$ is proved.

\smallskip

We give here a short proof of some statements in Step 3.

Observe that the topology of $K_{t_1,\dots,t_{n-1}}$ is essentially defined
as a topology of a finite-dimensional vector space over a standard field
$k_0\{\!\{t_1\}\!\}\dots\{\!\{t_{n-1}\}\!\}$. (It will be precisely so, if
we take $\{\pi_i e_j: 0\le i\le e-1,1\le j\le d\}$ as a basis of this vector space,
where $e$ is the absolute ramification index of $K$, and $\pi_{i+e}=p\pi_i$ for any $i$.)
This
enables one to reduce the statements to the case of a standard field $K$.

If $K$ is standard, then either $K=E((X))$ or $K=E\{\!\{X\}\!\}$, where
$E$ is of smaller dimension. Looking at expansions in $X$, it is easy to construct a limit
of any Cauchy sequence in $K$ and to prove the uniqueness of it. (In the case $K=E\{\!\{X\}\!\}$
one should use the alternative construction of topology in Remark 3 in {1.3.1}.)
This proves (1) in {1.3.2}.

To prove the sequential continuity of multiplication in the mixed characteristic case, let  $\alpha_i\to 0$ and $\beta_i\to 0$,
we shall show that $\alpha_i\beta_i\to 0$.

Since $\alpha_i\to 0,\beta_i\to 0$, one can easily see that there is $c\in\Bbb Z$ such that
$v_n(\alpha_i)\ge c$, $v_n(\beta_i)\ge c$ for $i\ge 1$.

By the above remark, we may assume that $K$ is standard, i.e., $K=E\{\!\{t\}\!\}$.
Fix an open subgroup $U$ in $K$; we have $P(d)\subset U$ for some integer $d$.
One can assume that $U=\Cal V_{\{V_i\}}$, $V_i$ are open subgroups in
$E$. Then there is $m_0$ such that $P_{E}(d-c) \subset V_m$ for $m>m_0$.
Let
$$
\alpha_i=\sum_{-\infty}^\infty a_i^{(r)} t^r, \quad
\beta_i=\sum_{-\infty}^\infty b_i^{(l)} t^l, \qquad
a_i^{(r)}, b_i^{(l)}\in E.
$$
Notice that one can find an $r_0$ such that
$a_i^{(r)}\in P_{E}(d-c)$ for $r<r_0$ and all $i$.
Indeed, if this were
not so, one could choose a sequence $r_1>r_2>\dots$ such
that $a_{i_j}^{(r_j)} \notin P_E(d-c)$ for some
$i_j$. It is easy to construct
a neighbourhood of zero $V'_{r_j}$ in $E$ such that
$P_E(d-c) \subset V'_{r_j}$,
$a_{i_j}^{(r_j)} \notin V_{r_j}$.
Now put $V'_r=E$ when $r$ is distinct from any of
$r_j$, and $U'=\Cal V_{\{V'_r\}}$. Then $a_{i_j} \notin U'$,
$j=1,2,\dots$ The set $\{i_j\}$ is obviously infinite, which 
contradicts  the condition $\alpha_i\to 0$.

Similarly, $b_i^{(l)} \in P_E(d-c)$ for $l<l_0$ and
all $i$. Therefore,
$$
\alpha_i \beta_i \equiv
\sum_{r=r_0}^{m_0} a_i^{(r)} t^r \cdot
\sum_{l=l_0}^{m_0} b_i^{(l)} t^l \mod U,
$$
and the condition $a_i^{(r)} b_i^{(l)} \to 0$ for all $r$
and $l$ immediately implies $\alpha_i\beta_i\to 0$.

\HHH 1.3.4. Expansion into power series

Let $n=2$. Then in characteristic $p$ we have
$\Bbb F_{q}\left( \left( X\right) \right) \left( \left( Y\right)
\right) =\{\sum \theta _{ij}X^{j}Y^{i}\}$,
\noindent where $\theta _{ij}$ are elements of $\Bbb F%
_{q}$ such that for some $i_0$ we have $\theta _{ij}=0$ for $i\le i_0$  and
for every $i$ there is $j(i)$ such that
$\theta _{ij}=0$ for $j \le j(i)$.

On the other hand, the definition of the topology implies  that
for every neighbourhood of zero $U$ there exists $i_0$ and  for every $i<i_0$
there exists $j(i)$ such that
$\theta X^{j}Y^{i}\in U$ whenever either $i\ge i_0$ or $i<i_0, j\ge j(i)$.

So every formal power series has only finitely many terms $\theta X^{j}Y^{i}$
outside $U$. Therefore, it is in fact a {\it convergent} power series in
the just defined topology.

\df Definition

  $\Omega \subset  \Bbb Z^{n}$ is called
admissible if for every $1 \le l \le n$ and every $j_{l+1},\dots ,j_{n}$
there is $i=i(j_{l+1},\dots ,j_{n})\in  \Bbb Z$ such that
$$
\left( i_{1},\dots ,i_{n}\right) \in \Omega ,\;i_{l+1}=j_{l+1},\dots
,i_{n}=j_{n}\Rightarrow i_{l} \ge i.
$$
\enddf

 \th Theorem

 Let $t_{1},\dots ,t_{n}$ be a system of local parameters
of $K$.  Let $s$ be a section of the residue map
$ {\CO} _{K}\rightarrow  {\CO }%
_{K}/ {\CM}_{K}$ such that  $s(0)=0$.  Let $\Omega $ be an admissible subset
of $\Bbb Z^{n}$.  Then the series
$$
\sum_{\left( i_{1},\dots ,i_{n}\right) \in \Omega }b_{i_{1},\dots
,i_{n}}t_{1}^{i_{1}}\dots t_{n}^{i_{n}}\;\text{converges\ }\quad (b_{i_{1},\dots
,i_{n}}\in s( {\CO} _{K}/ {\CM}_{K}))
$$
and every element of $K$ can be uniquely written this way.
\endth

\rk Remark

In this statement it is essential that the last residue field is finite.
In a more general setting, one should take a ``good enough'' section. For
example, for $K=k\left\{\! \left\{ T_{1}\right\}\! \right\} \dots\left\{\! \left\{ T_{m}\right\}\!
\right\} \allowmathbreak( ( T_{m+2}) ) \dots((T_{n}))$, where $k$ is a finite extension of the fraction field of $W(K_0)$ and $K_0$ is perfect of prime characteristic,
one may take the Teichm\"uller section $K_0\to K_{m+1}= k\left\{\! \left\{ T_{1}\right\}\! \right\} \dots\left\{\! \left\{ T_{m}\right\}\!\right\}$ composed with the obvious embedding $K_{m+1}\hookrightarrow K$.
\endrk

\pf Proof

We have
$$
\sum_{\left( i_{1},\dots ,i_{n}\right) \in \Omega }b_{i_{1},\dots,i_{n}}t_{1}^{i_{1}}\dots t_{n}^{i_{n}}
=
\sum_{b\in s( {\CO} _{K}/ {\CM}_{K})}
\bigl(b\cdot\sum_{\left( i_{1},\dots ,i_{n}\right) \in \Omega_b}t_{1}^{i_{1}}\dots t_{n}^{i_{n}}\bigr),
$$
where $\Omega_b=\{( i_{1},\dots ,i_{n}) \in \Omega : b_{i_{1},\dots,i_{n}}=b\}$.
In view of the property (4), it is sufficient to show that the inner sums converge.
Equivalently, one has to show that given a neighbourhood of zero $U$ in $K$,
for almost all $\left( i_{1},\dots ,i_{n}\right) \in \Omega$
 we have  $t_{1}^{i_{1}}\dots t_{n}^{i_{n}}\in U$.
This follows easily by induction on $n$ if we observe that
$t_{1}^{i_{1}}\dots t_{n-1}^{i_{n-1}} = h_{t_1,\dots,t_{n-1}}(\Overline{t_1}^{i_1}\dots\Overline{t_{n-1}}^{i_{n-1}})$.

To prove the second statement, apply induction on $n$ once again.
Let $r=v_n(\alpha)$, where $\alpha$ is a given element of $K$. Then
by the induction hypothesis
$$
\Overline{t_n^{-r}\alpha}=
\sum_{\left( i_{1},\dots ,i_{n-1}\right) \in \Omega_r }\Overline{b_{i_{1},\dots,i_{n}}}\bigl(\Overline{t_{1}}\bigr)^{i_{1}}\dots \bigl(\Overline{t_{n-1}}\bigr)^{i_{n-1}},
$$
where $\Omega_r\subset\Bbb Z^{n-1}$ is a certain admissible set. Hence
$$
\alpha=\sum_{\left( i_{1},\dots ,i_{n-1}\right) \in \Omega_r }
b_{i_{1},\dots,i_{n}}t_{1}^{i_{1}}\dots t_{n-1}^{i_{n-1}}t_n^r+\alpha',
$$
where $v_n(\alpha')>r$. Continuing this way, we obtain the desired expansion into a sum over
the admissible set $\Omega=(\Omega_r\times\{r\}) \cup (\Omega_{r+1}\times\{r+1\})\cup\dots$

The uniqueness follows from the continuity of the residue map $\Cal O_K\to K_{n-1}$.
\qed
\endpf

\HH 1.4. Topology on $K^*$

\HHH 1.4.1. 2-dimensional case, $\chr(k_K)=p$

\phantom{}\par\smallskip

Let $A$ be the last residue field $K_0$ if $\chr (K)=p$, and let $A=W(K_0)$
 if $\chr (K)=0$. Then $A$ is canonically embedded into $\Cal O_K$, and it
is in fact the subring generated by the set $\Cal R$.

For a 2-dimensional local field $K$  with a system of local parameters $t_2, t_1$ define a  base of neighbourhoods of 1 as
the set of all $1+t_2^i\Cal O_K+t_1^j A[[t_1,t_2]]$, $i\ge1$, $j\ge1$.
Then every element $\alpha\in K^*$ can be expanded as
a convergent (with respect to the just defined  topology)
product
$$\alpha=t_2^{a_2}t_1^{a_1}\theta\prod (1+\theta_{ij}t_2^i t_1^j)$$
with $\theta\in \Cal R^*, \theta_{ij}\in \Cal R, a_1,a_2\in\Bbb Z$.
The set $S=\{(j,i): \theta_{ij}\not=0\}$ is admissible.

\HHH 1.4.2

In the general case,  following Parshin's approach in characteristic $p$
\cite{P}, we define the topology $\tau$ on $K^*$ as follows.

\df Definition

If $\chr (K_{n-1})=p$, then define the topology $\tau$ on
$$K^*\simeq
V_K\times \langle t_1\rangle\times \dots\times \langle t_n\rangle \times \Cal R^*$$
as the product of the induced from $K$ topology on the group of principal units $V_K$ and  the discrete topology
on $\langle t_1\rangle\times \dots\times \langle t_n\rangle \times \Cal R^*$.

If $\chr (K)=\chr (K_{m+1})=0$, $\chr (K_{m})=p$, where $m\le n-2$, then
we have a canonical exact sequence
$$
1 @>>> 1+P_K\bigl(1,\underset{n-m-2}\to{\underbrace{0,\dots,0}}\,\bigr) @>>>
O _{K}^{*} @>>> O_{K_{m+1}}^* @>>> 1.
$$
Define
the  topology $\tau$  on
$K^*\simeq O _{K}^{*}\times \langle t_1\rangle\times \dots\times \langle t_n\rangle $ as
the product of the discrete topology on
$\langle t_{1}\rangle\times \dots\times \langle t_n\rangle$
and the inverse image of the topology $\tau$  on $O_{K_{m+1}}^*$.
Then the intersection of all neighbourhoods of 1
is equal to $1+P_K\bigl(1,\underset{n-m-2}\to{\underbrace{0,\dots,0}}\,\bigr)$
which is a uniquely divisible group.
\enddf

\rk Remarks

1. Observe that $K_{m+1}$ is a mixed characteristic field and therefore its
topology is well defined. Thus, the topology $\tau$ is well defined in all cases.

2. A base of neighbourhoods of 1 in $V_K$ is formed by the sets
$$
h(U_0)+h(U_1)t_n+...+h(U_{c-1})t_n^{c-1}+P_K(c),
$$
where $c\ge1$, $U_0$ is a neighbourhood of 1 in $V_{k_K}$, $U_1,\dots,U_{c-1}$ are neighbourhoods of zero in $k_K$,
$h$ is the canonical lifting associated with some local parameters, $t_n$ is the last local parameter of $K$.
In particular,  in the two-dimensional case $\tau$ coincides with the topology of {1.4.1}.
\endrk

\rk Properties

\Roster

\Item{(1)}
Each Cauchy sequence with respect to the topology $\tau$
converges in $K^*$.

\Item{(2)} Multiplication in $K^*$ is sequentially continuous.

\Item{(3)}  If $n\le 2$, then the multiplicative group $K^*$ is a topological group  and it has a countable base of open subgroups.
$K^*$ is not a topological group with respect to $\tau$
if $m\ge 3$.
\endRoster
\endrk

\pf Proof

(1) and (2) follow immediately from the corresponding properties of the topology
defined in subsection {1.3}. In the 2-dimensional case (3) is obvious from the
description given in 1.4.1. Next, let $m\ge3$, and let $U$ be an arbitrary neighbourhood
of 1. We may assume that $n=m$ and $U\subset V_K$. From the definition of the topology on $V_K$
we see that $U\supset 1+h(U_1)t_n+h(U_2)t_n^2$, where $U_1,U_2$ are neighbourhoods of $0$ in $k_K$,
$t_n$ a prime element in $K$, and $h$ the canonical lifting corresponding to some
choice of local parameters. Therefore,
$$\aligned
UU+P(4)&\supset (1+h(U_1)t_n)(1+h(U_2)t_n^2)+P(4)\\
&=\{1+h(a)t_n+h(b)t_n^2+h(ab)t_n^3: a\in U_1,b\in U_2\}+P(4).
\endaligned
$$
(Indeed, $h(a)h(b)-h(ab)\in P(1)$.)
Since $U_1U_2=k_K$ (see property (3) in {1.3.2}), it is clear that $UU$ cannot lie
in a neighbourhood of 1 in $V_K$ of the form
$1+h(k_K)t_n+h(k_K)t_n^2+h(U')t_n^3+P(4)$, where $U'\ne k_K$ is a
neighbourhood of 0 in $k_K$. Thus, $K^*$ is not a topological group.
\qed\endpf

\rk Remarks

1. From the point of view of class field theory and the existence theorem
one needs a stronger topology on $K^*$ than the topology $\tau$
 (in order to have more open subgroups).
For example,
for $n\ge3$ each open subgroup $A$ in $K^*$ with respect
to the topology $\tau$  possesses the property:
$1+t_n^2\Cal O_K\subset (1+t_n^3\Cal O_K)A$.

A topology $\lambda_*$ which is the sequential saturation of $\tau$
is introduced in subsection~6.2;
it has the same set of convergence sequences as $\tau$
but more open subgroups.
For example \cite{F1}, the subgroup in $1+t_n\Cal O_K$ topologically generated
by $1+\theta t_n^{i_n}\dots t_1^{i_1}$ with $(i_1,\dots,i_n)\not=
(0,0,\dots,1,2)$, $i_n\ge1$
(i.e., the sequential closure of the subgroup generated by these elements)
is open in $\lambda_*$ and
does not satisfy the above-mentioned property.

One can even introduce a topology on $K^*$ which has the same set of convergence sequences as $\tau$ and with respect to which $K^*$ is a topological group,
see \cite{F2}.

2. For another approach to define open subgroups of $K^*$
see the paper of K.~Kato  in this volume.
\endrk

\HHH 1.4.3. Expansion into convergent products

To simplify  the following statements we assume here $\chr k_K=p$.
Let $B$ be a fixed set of representatives of non-zero elements of the last residue field in $K$.

\th Lemma

Let $\{\alpha_i: i\in I\}$ be a subset of $V_K$ such that
$$
\alpha_i=1+\sum_{{\bold r}\in\Omega_i}b_{\bold r}^{(i)}t_1^{r_1}\dots t_n^{r_n},\tag{$*$}
$$
where $b\in B$, and $\Omega_i\subset\Bbb Z^n_+$ are admissible sets satisfying the following two conditions:

{\rm (i)} $\Omega=\bigcup_{i\in I} \Omega_i$ is an admissible set;

{\rm (ii)} $\bigcap_{j\in J}\Omega_j=\emptyset$, where $J$ is any infinite subset of $I$.

Then $\prod_{i\in I}\alpha_i$ converges.
\endth

\pf Proof

Fix a neighbourhood of 1 in $V_K$; by definition it is of the form $(1+U)\cap V_K$, where
$U$ is a neighbourhood of 0 in $K$.
Consider various finite products of $b_{\bold r}^{(i)}t_1^{r_1}\dots t_n^{r_n}$ which occur in \thetag{$*$}.
It is sufficient to show that almost all such products belong to $U$.

Any product under consideration has the form
$$
\gamma=b_1^{k_1}\dots b_s^{k_s}t_1^{l_1}\dots t_n^{l_n}\tag{$**$}
$$
with $l_n>0$, where $B=\{b_1,\dots,b_s\}$.
We prove by induction on $j$ the following claim: for $0\le j\le n$
and fixed $l_{j+1},\dots, l_n$ the element $\gamma$ almost always lies
in $U$
(in  case $j=n$ we obtain the original
claim).
Let $$
\hat\Omega=\{{\bold r}_1+\dots+{\bold r}_t :
t\ge1,{\bold r}_1,\dots,  {\bold r}_t\in\Omega\}.
$$
It is easy to see that $\hat\Omega$ is an admissible set and any element  of $\hat\Omega$
can be written as a sum of elements of $\Omega$ in finitely many ways only.
This fact and  condition (ii) imply that any particular $n$-tuple $(l_1,\dots,l_n)$
can occur at the right hand side of \thetag{$**$} only finitely many times. This proves the
base of induction ($j=0$).

For $j>0$, we see that $l_j$ is bounded from below since $(l_1,\dots,l_n)\in\hat\Omega$
and $l_{j+1},\dots,l_n$ are fixed.
On the other hand, $\gamma\in U$ for sufficiently large $l_j$
and arbitrary $k_1,\dots,k_s,l_1,\dots,l_{j-1}$ in view of \cite{MZ1, Prop. 1.4}
applied to the neighbourhood of zero $t_{j+1}^{-l_{j+1}}\dots t_n^{-l_n}U$ in $K$.
Therefore, we have to consider only
a finite range of values $c\le l_j\le c'$. For any $l_j$ in this range the induction
hypothesis is applicable.
\qed
\endpf

\th Theorem

For any ${\bold r}\in\Bbb Z^n_+$ and any $b\in B$ fix an element
$$
a_{{\bold r},b}=\sum_{{\bold s}\in\Omega_{{\bold r},b}}
b_{{\bold s}}^{{\bold r},b}t_1^{s_1}\dots t_n^{s_n},
$$
such that $b_{\bold r}^{{\bold r},b}=b$, and $b_{{\bold s}}^{{\bold r},b}=0$ for $\bold s<\bold r$.
Suppose that the admissible sets 
$$\{\Omega_{\bold r,b}: \bold r\in\Omega_*,b\in B\}$$
satisfy conditions {\rm(i)} and {\rm(ii)}
of the Lemma for any given admissible set $\Omega_*$.

1. Every element $a\in K$ can be uniquely expanded into a convergent
series
$$
a=\sum_{\bold r\in \Omega_a}a_{\bold r,b_{\bold r}},
$$
where $b_{\bold r}\in B$, $\Omega_a\subset\Bbb Z_n$ is an admissible set.

2. Every element $\alpha\in K^*$ can be uniquely expanded into a convergent product:
$$
\alpha =t_n^{a_{n}
}\dots t_1^{a_{1}}
b_0 \prod_{\bold r\in\Omega_\alpha}\bigl(1+a_{\bold r,b_{\bold r}}\bigr),
$$
where $b_0\in B$, $b_{\bold r}\in B$, $\Omega_\alpha\subset\Bbb Z_n^+$ is an admissible set.

\endth

\pf Proof

The additive part of the theorem is \cite{MZ2, Theorem 1}. The proof
of it is parallel to that of Theorem 1.3.4.

To prove the multiplicative part, we apply induction on $n$.
This reduces the statement to the case $\alpha\in 1+P(1)$. Here one
can construct an expansion and prove its uniqueness applying the
additive part of the theorem to the residue of $t_n^{-v_n(\alpha-1)}(\alpha-1)$ in $k_K$.
The convergence of all series which appear in this process follows from the above Lemma.
For details, see  \cite{MZ2, Theorem 2}.
\qed\endpf

\rk Remarks

1. Conditions (i) and (ii) in the Lemma are essential. Indeed, the infinite
products $\prod\limits_{i=1}^\infty(1+t_1^i+t_1^{-i}t_2)$ and
$\prod\limits_{i=1}^\infty(1+t_1^i+t_2)$ do not converge.
This means that the statements of Theorems 2.1 and 2.2 in \cite{MZ1} have to be corrected and
conditions (i) and (ii) for elements $\varepsilon_{\bold r,\theta}$ $(\bold  r\in\Omega_*)$ should be added.

2. If the last residue field is not finite, the statements are still true if the
system of representatives $B$ is not too pathological. For example, the system of
Teichm\"uller representatives is always suitable. The above proof works with the  only ammendment: instead
of Prop.\ 1.4 of \cite{MZ1} we apply the definition of topology directly.
\endrk

\th Corollary

If $\chr (K_{n-1})=p$, then every element $\alpha\in K^*$ can be expanded into a convergent product:
$$
\alpha =t_n^{a_{n}
}\dots t_1^{a_{1}}
\theta \prod (1+\theta_{i_{1},\dots
,i_{n}}t_1^{i_{1}}\dots
t_n^{i_{n}}), \quad \theta\in \Cal R^*,
\quad \theta_{i_{1},\dots,i_{n}}\in \Cal R,
\tag{$***$}
$$
with
$\{(i_1,\dots, i_n): \theta_{i_{1},\dots
,i_{n}}\not=0\}$ being an admissible set.  Any series \thetag{$***$} converges.
\endth

\Bib References

\rf{F1} I. Fesenko,
Abelian extensions of complete discrete valuation fields,
Number Theory Paris 1993/94,
Cambridge Univ. Press,  1996,
 47--74.

\rf{F2}
I.  Fesenko,
Sequential topologies and quotients of the Milnor $K$-groups of higher local
fields,
preprint, 
 www.maths.nott.ac.uk/personal/ibf/stqk.ps

\rf{FV}
 I.  Fesenko and S.  Vostokov,
Local Fields and Their Extensions,
AMS,
Providence,
1993.

\rf{K1}
K. Kato, A generalization of local class field theory
by using $K$-groups I,
J. Fac. Sci. Univ. Tokyo Sec. IA 26 No.2, 1979, 303--376.

\rf{K2}
K. Kato,
Existence Theorem for higher local class field theory,
this volume.

\rf{KZ}
 M. V. Koroteev and I. B. Zhukov,
 Elimination of wild ramification,
 Algebra i Analiz
11 (1999), no. 6, 153--177.

\rf{MZ1}
  A. I. Madunts and I. B. Zhukov,
 Multidimensional complete fields: topology
and other basic constructions,
 Trudy S.-Peterb. Mat. Obshch.
(1995);
English translation in
 Amer. Math. Soc. Transl. (Ser. 2)
165 (1995), 1--34.

\rf{MZ2}
A. I. Madunts and I. B. Zhukov,
Additive and multiplicative expansions
in multidimensional local fields,
Zap. Nauchn. Sem. POMI (to appear).

\rf{P}
 A. N. Parshin,
 Local class field theory,
 Trudy Mat. Inst. Steklov (1984);
English translation in Proc. Steklov Inst. Math. 165 (1985),  no. 3, 157--185.

\rf{Z}
 I. B. Zhukov,
 Structure theorems for complete fields,
 Trudy S.-Peterb. Mat. Obsch. (1995);
English  translation   in Amer. Math. Soc. Transl. (Ser. 2)
165 (1995), 175--192.

\endBib

\Coordinates

Department of Mathematics and Mechanics \
St. Petersburg University

Bibliotechnaya pl. 2, \
Staryj Petergof

198904 St. Petersburg \ Russia

E-mail: igor\@zhukov.pdmi.ras.ru

\endCoordinates

\vfill
\pagebreak
\end

%% file: m3-macs.tex
\expandafter\ifx\csname mthreemacsloaded\endcsname\relax\else \fi

\magnification1100
\input amstex


 \catcode`\@=11
 \let\wlog@ld\wlog
 \def\wlog#1{\relax}

 \newif\ifIN@
 \def\m@rker{\m@@rker}
 \def\IN@{\expandafter\INN@\expandafter}
 \long\def\INN@0#1@#2@{\long\def\NI@##1#1##2##3\ENDNI@
    {\ifx\m@rker##2\IN@false\else\IN@true\fi}%
     \expandafter\NI@#2@@#1\m@rker\ENDNI@}
  \newtoks\Initialtoks@  \newtoks\Terminaltoks@
  \def\SPLIT@{\expandafter\SPLITT@\expandafter}
  \def\SPLITT@0#1@#2@{\def\TTILPS@##1#1##2@{%
     \Initialtoks@{##1}\Terminaltoks@{##2}}\expandafter\TTILPS@#2@}
  \newtoks\Trimtoks@

 \def\ForeTrim@{\expandafter\ForeTrim@@\expandafter}
 \def\ForePrim@0 #1@{\Trimtoks@{#1}}
 \def\ForeTrim@@0#1@{\IN@0\m@rker. @\m@rker.#1@%
     \ifIN@\ForePrim@0#1@%
     \else\Trimtoks@\expandafter{#1}\fi}
 
  \def\Trim@0#1@{%
      \ForeTrim@0#1@%
      \IN@0 @\the\Trimtoks@ @%
        \ifIN@
             \SPLIT@0 @\the\Trimtoks@ @\Trimtoks@\Initialtoks@
             \IN@0\the\Terminaltoks@ @ @%
                 \ifIN@
                 \else \Trimtoks@ {FigNameWithSpace}%
                 \fi
        \fi
      }

  \font\titlebold=cmbx12 scaled 1200
  \font\twelvebold=cmbx12
  \font\tenbold=cmbx10
  \font\ninebold=cmbx9
  \font\sevenbold=cmbx7
  \font\fivebold=cmbx5

  \input amssym.def \input amssym
     \font\titlemsa=msam10 at 14.4pt
     \font\titlemsb=msbm10 at 14.4pt
     \font\titleeufm=eufm10 at 14.4pt
     \font\twelvemsa=msam10 scaled 1200
     \font\twelvemsb=msbm10 scaled 1200
     \font\twelveeufm=eufm10 scaled 1200
     \font\ninemsa=msam9
     \font\ninemsb=msbm9
     \font\nineeufm=eufm9

   \ifx\cyrfam\undefined
   \else
     \immediate\write16{}%
     \message{ !!! cyr fonts already defined. !!! }
     \message{ --- edit out superfluous font defs? }
   \fi
   \newfam\cyrfam
       \font\titlecyr=wncyr10 scaled 1440 
       \font\twelvecyr=wncyr10 scaled 1200
       \font\tencyr=wncyr10
       \font\ninecyr=wncyr9
       \font\sevencyr=wncyr7
       \font\sixcyr=wncyr6

   \newfam\eusmfam
       \font\titleeusm=eusm10 scaled 1440
       \font\twelveeusm=eusm10 scaled 1200
       \font\teneusm=eusm10
       \font\nineeusm=eusm9
       \font\seveneusm=eusm7
       
       \font\fiveeusm=eusm5

\let\Cal\cal

    \font\ninemrm=cmr9 
    \font\ninei=cmmi9
    \font\ninesy=cmsy9 
    \skewchar\ninei='177
    \skewchar\ninesy='60

  \font\twelvemrm=cmr10 at 12pt 
  \font\twelvei=cmmi10 at 12pt
  \font\twelvesy=cmsy10 at 12pt

  \font\titlemrm=cmr10 at 14.4pt 
  \font\titlei=cmmi10 at 14.4pt
  \font\titlesy=cmsy10 at 14.4pt


  \def\Smallfonts{\ninepoint}

  \def\Hfont{\titlepoint\bf}
  \def\Authorfont{\twelvepoint\it}
  \def\HHfont{\twelvepoint\bf}
  \def\HHHfont{\bf}
  \def\Bibfont{\tenbf}
  \def\Coordfont{\nineit }

  \def \thfont {\bf }
  \def \pffont {\it\itSpacing }
  \def \rkfont {\bf }
  \def \dffont {\bf }
  \def \egfont {\bf }

 \def\ninepoint{%
  \def\rm{\fam0\ninerm}%
    \textfont0=\ninemrm  \scriptfont0=\sevenrm  \scriptscriptfont0=\fiverm
    \textfont1=\ninei    \scriptfont1=\seveni   \scriptscriptfont1=\fivei
  \def\mit{\fam1\ninei}%
  \def\oldstyle{\fam1\ninei}%
    \textfont2=\ninesy   \scriptfont2=\sevensy  \scriptscriptfont2=\fivesy
    \textfont3=\tenex    \scriptfont3=\tenex    \scriptscriptfont3=\tenex
  \def\it{\fam\itfam\nineit}%
    \textfont\itfam=\nineit
  \def\bf{\ifmmode\fam\bffam\else\ninebf\fi}%
    \textfont\bffam=\ninebold 
    \scriptfont\bffam=\sevenbold 
    \scriptscriptfont\bffam=\fivebold%
  \def\msa{\fam\msafam\ninemsa}%
    \textfont\msafam=\ninemsa 
    \scriptfont\msafam=\sevenmsa
    \scriptscriptfont\msafam=\fivemsa%
  \def\msb{\fam\msbfam\ninemsb}%
    \textfont\msbfam=\ninemsb%
    \scriptfont\msbfam=\sevenmsb%
    \scriptscriptfont\msbfam=\fivemsb%
  \def\eufm{\fam\eufmfam\nineeufm}%
    \textfont\eufmfam=\nineeufm
    \scriptfont\eufmfam=\seveneufm
    \scriptscriptfont\eufmfam=\fiveeufm
   \def\eusm{\fam\eusmfam\nineeusm}%
     \textfont\eusmfam=\nineeusm
     \scriptfont\eusmfam=\seveneusm
     \scriptscriptfont\eusmfam=\fiveeusm
   \def\cyr{\fam\cyrfam\ninecyr}%
     \textfont\cyrfam=\ninecyr
     \scriptfont\cyrfam=\sevencyr
     \scriptscriptfont\cyrfam=\sixcyr
  \setbox\strutbox=\hbox{\vrule
      height7pt depth3pt width0pt}%
   \baselineskip=10.8pt\rm}

 \let\eightpoint\ninepoint 

 \def\tenpoint{%
  \def\rm{\fam0\tenrm}%
    \textfont0=\tenmrm \scriptfont0=\sevenrm \scriptscriptfont0=\fiverm%
  \def\mit{\fam1\teni}%
  \def\oldstyle{\fam1\teni}%
    \textfont1=\teni   \scriptfont1=\seveni  \scriptscriptfont1=\fivei%
    \textfont2=\tensy  \scriptfont2=\sevensy \scriptscriptfont2=\fivesy%
    \textfont3=\tenex  \scriptfont3=\tenex   \scriptscriptfont3=\tenex%
  \def\it{\fam\itfam\tenit}%
    \textfont\itfam=\tenit%
  \def\bf{\ifmmode\fam\bffam\else\tenbf\fi}%
    \textfont\bffam=\tenbold
    \scriptfont\bffam=\sevenbold%
    \scriptscriptfont\bffam=\fivebold%
  \def\msa{\fam\msafam\tenmsa}%
    \textfont\msafam=\tenmsa%
    \scriptfont\msafam=\sevenmsa%
    \scriptscriptfont\msafam=\fivemsa%
  \def\msb{\fam\msbfam\tenmsb}%
    \textfont\msbfam=\tenmsb%
    \scriptfont\msbfam=\sevenmsb%
    \scriptscriptfont\msbfam=\fivemsb%
  \def\eufm{\fam\eufmfam\teneufm}%
   \textfont\eufmfam=\teneufm
   \scriptfont\eufmfam=\seveneufm
   \scriptscriptfont\eufmfam=\fiveeufm
   \def\eusm{\fam\eusmfam\teneusm}%
    \textfont\eusmfam=\teneusm
    \scriptfont\eusmfam=\seveneusm
    \scriptscriptfont\eusmfam=\fiveeusm
   \def\cyr{\fam\cyrfam\tencyr}%
    \textfont\cyrfam=\tencyr
    \scriptfont\cyrfam=\sevencyr
    \scriptscriptfont\cyrfam=\sixcyr
  \setbox\strutbox=\hbox{\vrule %
      height8.5pt depth3.5ptwidth0pt}%
  \baselineskip=\StdBaselineskip\rm}

 \def\twelvepoint{%
  \def\rm{\fam0\twelverm}%
    \textfont0=\twelvemrm \scriptfont0=\tenmrm \scriptscriptfont0=\sevenrm
    \textfont1=\twelvei   \scriptfont1=\teni   \scriptscriptfont1=\seveni
  \def\mit{\fam1\twelvei}%
  \def\oldstyle{\fam1\twelvei}%
    \textfont2=\twelvesy  \scriptfont2=\tensy  \scriptscriptfont2=\sevensy
    \textfont3=\tenex  \scriptfont3=\tenex  \scriptscriptfont3=\tenex
  \def\it{\fam\itfam\twelveit}%
    \textfont\itfam=\twelveit
  \def\bf{\ifmmode\fam\bffam\else\twelvebf\fi}%
    \textfont\bffam=\twelvebold
    \scriptfont\bffam=\tenbold%
    \scriptscriptfont\bffam=\sevenbold%
  \def\msa{\fam\msafam\twelvemsa}%
    \textfont\msafam=\twelvemsa%
    \scriptfont\msafam=\tenmsa%
    \scriptscriptfont\msafam=\sevenmsa%
  \def\msb{\fam\msbfam\twelvemsb}%
    \textfont\msbfam=\twelvemsb%
    \scriptfont\msbfam=\tenmsb%
    \scriptscriptfont\msbfam=\sevenmsb%
  \def\eufm{\fam\eufmfam\twelveeufm}%
   \textfont\eufmfam=\twelveeufm
   \scriptfont\eufmfam=\teneufm
   \scriptscriptfont\eufmfam=\seveneufm
   \def\eusm{\fam\eusmfam\twelveeusm}%
    \textfont\eusmfam=\twelveeusm
    \scriptfont\eusmfam=\teneusm
    \scriptscriptfont\eusmfam=\seveneusm
   \def\cyr{\fam\cyrfam\tencyr}%
    \textfont\cyrfam=\twelvecyr
    \scriptfont\cyrfam=\tencyr
    \scriptscriptfont\cyrfam=\sevencyr
  \setbox\strutbox=\hbox{\vrule
      height10.2pt depth4.55pt width0pt}%
  \baselineskip=14pt\rm}

 \def\titlepoint{%
    \textfont0=\titlemrm \scriptfont0=\twelvemrm \scriptscriptfont0=\tenmrm
    \textfont1=\titlei   \scriptfont1=\twelvei   \scriptscriptfont1=\teni
  \def\mit{\fam1\titlei}%
  \def\oldstyle{\fam1\titlei}%
    \textfont2=\titlesy  \scriptfont2=\twelvesy  \scriptscriptfont2=\tensy
    \textfont3=\tenex
    \scriptfont3=\tenex
    \scriptscriptfont3=\tenex
  \def\it{\fam\itfam\titleit}%
    \textfont\itfam=\titleit
  \def\bf{\ifmmode\fam\bffam\else\titlebf\fi}%
    \textfont\bffam=\titlebold
    \scriptfont\bffam=\twelvebold%
    \scriptscriptfont\bffam=\tenbold%
  \def\msa{\fam\msafam\titlemsa}%
    \textfont\msafam=\titlemsa%
    \scriptfont\msafam=\twelvemsa%
    \scriptscriptfont\msafam=\tenmsa%
  \def\msb{\fam\msbfam\titlemsb}%
    \textfont\msbfam=\titlemsb%
    \scriptfont\msbfam=\twelvemsb%
    \scriptscriptfont\msbfam=\tenmsb%
  \def\eufm{\fam\eufmfam\titleeufm}%
    \textfont\eufmfam=\titleeufm
    \scriptfont\eufmfam=\twelveeufm
    \scriptscriptfont\eufmfam=\teneufm
   \def\eusm{\fam\eusmfam\titleeusm}%
     \textfont\eusmfam=\titleeusm
     \scriptfont\eusmfam=\twelveeusm
     \scriptscriptfont\eusmfam=\teneusm
   \def\cyr{\fam\cyrfam\tencyr}%
    \textfont\cyrfam=\titlecyr
    \scriptfont\cyrfam=\twelvecyr
    \scriptscriptfont\cyrfam=\tencyr
  \setbox\strutbox=\hbox{\vrule
      height12.3pt depth5.54pt width0pt}%
  \baselineskip=16pt\rm}

\newbox\AuthorBox\newbox\TitleBox
\newbox\TFLinebox
\newbox\FLinebox
\newbox\HLinebox
\def\SetTFLinebox#1{\setbox\TFLinebox=\hbox{#1}}
\def\SetFLinebox#1{\setbox\FLinebox=\hbox{#1}}
\def\SetHLinebox#1{\setbox\HLinebox=\hbox{#1}}

 \def\SetAuthorHead#1{%
     \setbox\AuthorBox=\hbox{\ninepoint \it 
           \ignorespaces\frenchspacing#1\unskip}}
 \def\SetTitleHead#1{%
     \setbox\TitleBox=\hbox{\ninepoint \it
           \ignorespaces\frenchspacing#1\unskip}}

  \def\itSpacing{\relax}
  \def\itSpacingOff{\relax}


 \def\Hrule{\hrule width0pt height0pt}

  \newskip\ProcSkip \ProcSkip 8pt plus2pt minus2pt

 \newskip\LastSkip
 \def\SaveLastSkip{\LastSkip\lastskip}
 \def\RestoreLastSkip{\vskip-\LastSkip\vskip\LastSkip}

 \def\NoindentAfter{\everypar={\setbox0=\lastbox\everypar={}}}

 \long\def\H#1\par#2\par{\notenumber=0 \titlepagetrue%
    {
    \baselineskip=20pt
    \parindent=0pt\parskip=0pt\frenchspacing
    \leftskip=0pt plus .2\hsize minus .3\hsize
    \rightskip=0pt plus .2\hsize minus .3\hsize
 \def\\{\unskip\break}%
    \pretolerance=10000 \Hfont #1\unskip\break
     \vskip7pt\Hrule
\hfill \Authorfont #2\hfill\hfill\unskip}
    \vskip48pt plus 4pt minus 4pt
    \par\NoindentAfter\rm}

 \long\def\Hi#1\par#2\par{\notenumber=0 \titlepagetrue%
    {  \baselineskip=0pt  \parindent=0pt\parskip=0pt\frenchspacing
    \leftskip=0pt plus .2\hsize minus .3\hsize
    \rightskip=0pt plus .2\hsize minus .3\hsize
}
    \rm}


 \newdimen\PageRemainder
  \def\SetPageRemainder{
     \PageRemainder=\pagegoal
     \ifdim\PageRemainder=\maxdimen\PageRemainder=\vsize
     \else\advance\PageRemainder by -1\pagetotal\fi}

  \def\Rpt@{}\def\Rpt@@{}

  \long\def\HH#1\par{\par
  \SaveLastSkip\removelastskip\goodbreak
  \ifdim\LastSkip<30pt 
     \LastSkip 30pt
plus 3pt minus 2pt\fi
  \SetPageRemainder\advance\PageRemainder-\LastSkip
  \ifdim\PageRemainder<150pt
       \edef\Rpt@{remain = \the\PageRemainder\noexpand\\
                pagetotal=\the\pagetotal\noexpand\\
                           pagegoal=\the\pagegoal}%
          \fi
   \ifdim\PageRemainder<65pt 
       \ifdim\PageRemainder > 0pt
          \edef\Rpt@@{\noexpand\\
                      Had HH PageRemainder$<$\relax 65pt\noexpand\\
                      Hence forced break!}%
     \vskip 0pt plus .2\PageRemainder\eject 
    \fi\fi
    \vskip\LastSkip\Hrule 
    \pretolerance=10000\rightskip=0pt plus 3em
    \hangafter1 \hangindent=2.2em%
    \noindent
    \HHfont \unskip \Ednote{\Rpt@\Rpt@@}%
            \def\Rpt@{}\def\Rpt@@{}%
            \ignorespaces
            #1\par\rightskip=0pt\pretolerance=\StdPretolerance%
    \NoindentAfter
\tenpoint\rm%
     \medskip \vskip\ProcSkip}

  \long\def\HHH#1\par{\par%
  \SaveLastSkip\removelastskip\goodbreak
  \ifdim\LastSkip<\ProcSkip%
     \LastSkip\ProcSkip\fi
  \SetPageRemainder\advance\PageRemainder-\LastSkip
  \ifdim\PageRemainder<150pt
       \edef\Rpt@{remain = \the\PageRemainder\noexpand\\
                pagetotal=\the\pagetotal\noexpand\\
                           pagegoal=\the\pagegoal}%
       \fi
   \ifdim\PageRemainder<48pt  
        \ifdim\PageRemainder > 0pt
             \edef\Rpt@@{\noexpand\\
                      Had HHH PageRemainder$<$\relax48pt\noexpand\\
                      Hence forced break!}%
       \vskip 0pt plus .2\PageRemainder\eject 
      \fi\fi
   \vskip\LastSkip\par\noindent
   \HHHfont \unskip\Ednote{\Rpt@\Rpt@@}%
  \def\Rpt@{}\def\Rpt@@{}%
  \ignorespaces
   #1\unskip.\quad\rm\ignorespaces
   \ignorepars}

  \long\def\ignorepars#1\par{\def\Test{#1}%
     \ifx\Test\Empty\def\This{\ignorepars}%
        \else\def\This{\Test\par}\fi
           \This}
  \def\Empty{}

 \def\Abstract#1\par{\bgroup\Smallfonts\narrower\HHH #1\par}
 \def\endAbstract{\par\egroup}


 \def\ProcBreak{\par%
    \ifdim\lastskip<8pt%
    \removelastskip%
    \penalty-200\vskip\ProcSkip\fi}

 \def\th#1\par{\ProcBreak \noindent
   {\thfont\ignorespaces
    #1\unskip.}\it\itSpacing\kern.4em\ignorepars}

 \def\endth{\ProcBreak\rm\itSpacingOff }


 \def\pf#1\par{\ProcBreak %
    \noindent\pffont#1\unskip.\rm\itSpacingOff{\kern .7em}\ignorepars}

 \def\endpf{\medskip \ProcBreak } 

  \def\qedbox{\hbox{\vbox{
    \hrule width0.2cm height0.2pt
    \hbox to 0.2cm{\vrule height 0.2cm width 0.2pt
             \hfil\vrule height0.2cm width 0.2pt}
    \hrule width0.2cm height 0.2pt}\kern1pt}}

  \def\qed{\ifmmode\qedbox
    \else\unskip\ \hglue0mm\hfill\qedbox\ProcBreak\fi}

  \def \rk #1\par{\ProcBreak
     \noindent{\rkfont\ignorespaces #1\unskip.}%
     \rm\kern.6em\ignorepars}

  \def \endrk {\medskip\ProcBreak }

  \def \df #1\par{\ProcBreak
     \noindent{\dffont\unskip\ignorespaces #1\unskip.}%
     \rm\kern.6em\ignorepars}

  \def \enddf {\medskip\ProcBreak }

  \def \eg #1\par{\ProcBreak
     \noindent\egfont\unskip\ignorespaces #1\unskip.
     \rm\kern.6em\ignorepars}

  \def \endeg {\medskip\ProcBreak }

  \newdimen\Overhang

   \def\MaxTag@#1#2#3#4#5{\setbox0=\hbox{#4\ignorespaces#2\unskip}%
     \dimen0=\wd0\advance\dimen0 by#3
     \ifdim\dimen0<#5\relax\dimen0=#5\fi
     \expandafter\edef\csname #1Hang\endcsname{\the\dimen0}}

 \def\MaxItemTag#1{\MaxTag@{Item}{#1}{.4em}{\ItemStyle}{\parindent}}%
 \def\MaxItemItemTag#1{%
        \MaxTag@{ItemItem}{#1}{.4em}{\ItemItemStyle}{\parindent}}
 \def\MaxNrTag#1{\MaxTag@{Nr}{#1}{.5em}{\NrStyle}{\parindent}}
 \def\MaxReferenceTag#1{%
        \MaxTag@{Reference}{[#1]}{.6em}{\ninerm}{\parindent}}
 \def\MaxFootTag#1{\MaxTag@{Foot}{#1}{.4em}{\ninerm}{\z@}}

  \def\SetOverhang@{\Overhang=.8\dimen0%
     \advance\Overhang by \wd0\relax
     \ifdim\Overhang>\hangindent\relax
       \advance\Overhang by .25\dimen0%
       \Ednote{Tag is pushing text.}\osumess{Tag is pushing text.}%
     \else\Overhang=\hangindent
     \fi}

   \def\Item#1{\par\noindent
      \hangafter1\hangindent=\ItemHang
      \setbox0=\hbox{\ItemStyle\ignorespaces#1\unskip}%
      \dimen0=.4em\SetOverhang@
      \rlap{\box0}\kern\Overhang\ignorespaces}

   \def\ItemItem#1{\par\noindent
      \hangafter1\hangindent=\ItemItemHang
      \setbox0=\hbox{\ItemItemStyle\ignorespaces#1\unskip}%
      \dimen0=.4em\SetOverhang@
      \advance\hangindent by \ItemHang
      \kern\ItemHang\rlap{\box0}%
      \kern\Overhang\ignorespaces}

  \def\Nr#1{\par\noindent\hangindent=\NrHang 
    \setbox0=\hbox{\NrStyle\ignorespaces#1\unskip}%
    \dimen0=.5em\SetOverhang@
    \rlap{\box0}\kern\Overhang
    \hangindent=\z@\ignorespaces}

   \newskip\Rosterskip\Rosterskip 1pt plus1pt 
   \def\Roster{\par\ifdim\lastskip<\Rosterskip\removelastskip\vskip\Rosterskip\fi
    \bgroup}
   \def\endRoster{\par\global\edef\LastSkip@{\the\lastskip}\removelastskip
       \egroup\penalty-50\LastSkip\LastSkip@\relax
       \ifdim\LastSkip<\Rosterskip\LastSkip\Rosterskip\fi
       \vskip\LastSkip}




 \def\cite#1{
    \def\nextiii@##1,##2\end@{{\frenchspacing\rm 
      \lBr\ignorespaces##1\unskip{\rm,~\ignorespaces##2}\rBr}}%
    \IN@0,@#1@%
    \ifIN@\def\next{\nextiii@#1\end@}\else
    \def\next{{\rm\lBr#1\rBr}}\fi\next}


   \def \Bib#1\par{%
       \par\removelastskip\SetPageRemainder
       \ifdim\PageRemainder < 97pt
        \ifdim\PageRemainder > 0pt
        \vfill\eject
       \fi\fi
    \ProcBreak \par\begingroup\parskip=0 pt%
    \goodbreak \vskip 15 pt plus 10 pt
    \noindent\null\hfill\Bibfont
      \ignorespaces #1\unskip\hfill\null\par 
    \frenchspacing \Smallfonts\rm
    \parskip=2.5 pt plus 1 pt minus.5pt%
    \nobreak\vskip 12pt plus 2pt minus2pt\nobreak
    \leftskip=0 pt \baselineskip=10.5pt}

 \def\ReferenceTagSlide{0em}
  \def\ReferenceTagGap{.5em}

  \def \rf#1{\par\noindent
     \hangafter1\hangindent=\ReferenceHang      
     \setbox0=\hbox{\ninerm[\ignorespaces#1\unskip]}%
     \dimen0=\ReferenceTagGap\SetOverhang@
     \rlap{\kern\ReferenceTagSlide\box0}%
     \kern\Overhang\ignorespaces}

  \def\ref#1\par#2\par#3\par#4\par{%
     \rf{#1}#2\unskip,\ #3\unskip,\
     #4\unskip.}

  \def\endBib{\par\endgroup\vskip 12pt minus 6pt }


  \long\def\Coordinates#1\endCoordinates{
 {\par\vskip4pt\def\\{\unskip, }\Coordfont\baselineskip10.5pt\noindent#1}}

 \def\pagecontents{
  \gdef\Pagetot@l{\pagetotal}
  \ifvoid\TRMargIns\else
    \rlap{\kern\hsize\kern10pt\vbox to 0pt{%
         \box\TRMargIns\vss}}\fi
  \ifvoid\topins\else\unvbox\topins\fi
   \dimen@=\dp\@cclv \unvbox\@cclv 
   \ifvoid\footins\else 
     \vskip\skip\footins
     \footnoterule
     \unvbox\footins\fi
   \ifr@ggedbottom \kern-\dimen@ \vfil \fi}


 \newcount\Ht 

 \def \Acc{\expandafter } 

 \def\swthat{\raise -1.1 ex\hbox{\sam$\widehat{}$}}
 \def\swttilde{\raise -1.2 ex\hbox{\sam$\widetilde{}$}}
 \def \overdot{{\raise .2 ex \hbox to 0pt {\hss\bf\smash{.}\hss}}}
 \def \overcircle{{\raise .1 ex \hbox to 0pt
    {\sam$\eightpoint\scriptstyle\hss\circ\hss$}}}

 \def \Mathaccent#1#2{{\sam 
  \setbox4=\hbox{$\vphantom{#2}$}
  \Ht=\ht4 
  \setbox5=\hbox{${#1}$}
  \setbox6=\hbox{${#2}$}
  \setbox7=\hbox to .5\wd6{}
  \copy7\kern .1\Ht \raise\Ht sp\hbox{\copy5}\kern-.1\Ht
  \copy7\llap{\box6}
  }}

  \def\SwtCheck #1{
        \ifmmode \check{#1}%
                \else \v {#1}%
                \fi}

 \def\barpartial {%
   \kern .17 em
    \overline {\kern -.17 em\partial\kern-.03 em}%
    \kern .03 em}

 
  \def\Overline#1{\setbox1=\hbox{\sam ${#1}$}%
      \ifdim \wd1 > 6pt
    \kern .11 em
    \overline {\kern -.11 em#1\kern-.14 em}
    \kern .14 em
  \else
    \kern .03 em
    \overline {\kern -.03 em#1\kern-.04 em}
    \kern .04 em
  \fi}

 \def\SOverline#1{\setbox1=\hbox{\sam ${#1}$}%
      \ifdim \wd1 > 7pt
    \kern .22 em
    \overline {\kern -.22 em#1\kern-.09 em}%
    \kern .09 em
  \else
    \kern .10 em
    \overline {\kern -.10 em#1\kern-.04 em}%
    \kern .04 em
  \fi}


 \def\Underline#1{\setbox1=\hbox{\sam ${#1}$}%
      \ifdim \wd1 > 6pt
    \kern .11 em
    \underline {\kern -.11 em#1\kern-.14 em}
    \kern .14 em
  \else
    \kern .03 em
    \underline {\kern -.03 em#1\kern-.04 em}
    \kern .04 em
  \fi}

 \def\SUnderline#1{\setbox1=\hbox{\sam ${#1}$}%
      \ifdim \wd1 > 7pt
    \kern .04 em
    \underline {\kern -.04 em#1\kern-.2 em}%
    \kern .2 em
  \else
    \kern .0 em
    \underline {\kern -.0 em#1\kern-.15 em}%
    \kern .15 em
  \fi}


 \def \Blackbox
   {\leavevmode\hskip .3pt \vbox
   {\hrule height 5pt\hbox{\hskip 4.5pt}}\hskip .5pt}

 \def \XX{\Blackbox\kern.5pt\Blackbox} 

  \def\.{.\kern1pt}

    \def\Hyphen{\edef\this{\the\hyphenchar\font}%
          \hyphenchar\font=-1\char\this\hyphenchar\font=\this}

 \ifx\undefined\text
  \def\text#1{\hbox{\rm #1}}\fi 



   \everymath{}  

  \def\PassMath@@{\aftergroup\AfterMath@} 

 \let\PassMath@\PassMath@@

 \def\AfterMath@{\futurelet\next\AfterMathMole@}

 \def\AfterMathMole@{
      \ifcat\next\space
          \def\this{}
      \else
      \ifcat\next\egroup %
        \def\this{\osumess{Handset mathsurround?? ---(see dollar brace)}}%
      \else
      \def\this{\AAfterMath@}
      \fi\fi
      \this}

 \def\hyphen@{-}
 \def\paren@{)}
 \def\apostr@{'}

 \def\MSC#1{\kern-.8\mathsurround#1\kern.8\mathsurround}

 \def\AAfterMath@#1{\def\Next{#1}
    \IN@0\Next @,.;:!?\relax @%
    \ifIN@\def\this{\MSC{\Next}}%
    \else
    \ifx\Next\hyphen@\def\this{\futurelet\next\AfterHyphen@}%
    \else
    \ifx\Next\paren@\def\this{#1}%
    \else 
    \ifx\Next\apostr@\def\this{#1}%
    \else \def\this{\osumess{Handset mathsurround??}%
                 #1}\fi\fi\fi\fi
    \this}

 \def\AfterHyphen@#1{\def\Next{#1}%
   \ifx\Next\hyphen@\def\this{--}\else
   \ifcat\next\space%
   \def\this{\kern-\mathsurround\kern.05em- \Next}\else
   \def\this{\kern-\mathsurround\kern.05em\Hyphen\Next}\fi\fi\this}

 \def\sam{\mathsurround=\z@\let\PassMath@\relax}  %
 \def\mas{\mathsurround=\StdMathsurround\let\PassMath@\PassMath@@}
 
 \def\Mas{\mathsurround=\StdMathsurround
                \everymath{\PassMath@}\let\PassMath@\PassMath@@}

 \def\m@th{\mathsurround=\z@\everymath{}}

 \def\m@@th{\mathsurround=\z@\everymath={}\let\m@th\relax}

\def\underbar#1{$\setbox\z@\hbox{#1}\dp\z@\z@
      \m@th \underline{\box\z@}$\relax}

\def\mathhexbox#1#2#3{\leavevmode
  \hbox{\m@@th$\m@th \mathchar"#1#2#3$}}

\def\dots{\relax\ifmmode\ldots\else$\m@th\ldots\,$\relax\fi}

\def\dotfill{\cleaders\hbox{\m@@th$\m@th \mkern1.5mu.\mkern1.5mu$}\hfill}
\def\rightarrowfill{$\m@th\mathord-\mkern-6mu%
  \cleaders\hbox{\m@@th$\mkern-2mu\mathord-\mkern-2mu$}\hfill
  \mkern-6mu\mathord\rightarrow$\relax}
\def\leftarrowfill{$\m@th\mathord\leftarrow\mkern-6mu%
  \cleaders\hbox{\m@@th$\mkern-2mu\mathord-\mkern-2mu$}\hfill
  \mkern-6mu\mathord-$\relax}

\def\downbracefill{$\m@th\braceld\leaders\vrule\hfill\braceru
  \bracelu\leaders\vrule\hfill\bracerd$\relax}
\def\upbracefill{$\m@th\bracelu\leaders\vrule\hfill\bracerd
  \braceld\leaders\vrule\hfill\braceru$\relax}

\def\angle{{\vbox{\m@@th\ialign{$\m@th\scriptstyle##$\crcr
      \not\mathrel{\mkern14mu}\crcr
      \noalign{\nointerlineskip}
      \mkern2.5mu\leaders\hrule height.34pt\hfill\mkern2.5mu\crcr}}}}

\def\big#1{{\m@@th\hbox{$\left#1\vbox to8.5\p@{}\right.\n@space$}}}
\def\Big#1{{\m@@th\hbox{$\left#1\vbox to11.5\p@{}\right.\n@space$}}}
\def\bigg#1{{\m@@th\hbox{$\left#1\vbox to14.5\p@{}\right.\n@space$}}}
\def\Bigg#1{{\m@@th\hbox{$\left#1\vbox to17.5\p@{}\right.\n@space$}}}
\def\n@space{\nulldelimiterspace\z@ \m@th}

\def\root#1\of{\setbox\rootbox\hbox{\m@@th$\m@th\scriptscriptstyle{#1}$}
  \mathpalette\r@@t}
\def\r@@t#1#2{\setbox\z@\hbox{\m@@th$\m@th#1\sqrt{#2}$\relax}
  \dimen@\ht\z@ \advance\dimen@-\dp\z@
  \mkern5mu\raise.6\dimen@\copy\rootbox \mkern-10mu \box\z@}

\def\mathph@nt#1#2{\setbox\z@\hbox{\m@@th$\m@th#1{#2}$}\finph@nt}

\def\mathsm@sh#1#2{\setbox\z@\hbox{\m@@th$\m@th#1{#2}$}\finsm@sh}

\def\@vereq#1#2{\lower.5\p@\vbox{\m@@th\baselineskip\z@skip\lineskip-.5\p@
    \ialign{$\m@th#1\hfil##\hfil$\crcr#2\crcr=\crcr}}}

\def\mathpalette#1#2{\sam\mathchoice{#1\displaystyle{#2}}%
  {#1\textstyle{#2}}{#1\scriptstyle{#2}}{#1\scriptscriptstyle{#2}}\mas}

\def\widehat#1{\setbox\z@\hbox{\sam$#1$}%
 \ifdim\wd\z@>\tw@ em\mathaccent"0\msbfam@5B{#1}%
 \else\mathaccent"0362{#1}\fi}
\def\widetilde#1{\setbox\z@\hbox{\sam$#1$}%
 \ifdim\wd\z@>\tw@ em\mathaccent"0\msbfam@5D{#1}%
 \else\mathaccent"0365{#1}\fi}

 \def\dots{\relax{}
  \ifmmode\def\thedots{\mdots@}\else\def\thedots{\tdots@}\fi %
  \thedots}

 \let\@ldeqno\eqno\let\@ldleqno\leqno
 \def\eqno{\everymath{}\@ldeqno} \def\leqno{\everymath{}\@ldleqno}

  \let\@ldeqalignno\eqalignno
  \def\eqalignno#1{\sam\@ldeqalignno{#1}\mas}
  \let\@ldeqalign\eqalign
  \def\eqalign#1{\sam\@ldeqalign{#1}\mas}

 \def\overrightarrow#1{\vbox{\m@th\ialign{##\crcr
      \rightarrowfill\crcr\noalign{\kern-\p@\nointerlineskip}
      $\hfil\displaystyle{#1}\hfil$\crcr}}}
 \def\overleftarrow#1{\vbox{\m@th\ialign{##\crcr
      \leftarrowfill\crcr\noalign{\kern-\p@\nointerlineskip}
      $\hfil\displaystyle{#1}\hfil$\crcr}}}
 \def\overbrace#1{\mathop{\vbox{\m@th\ialign{##\crcr\noalign{\kern3\p@}
      \downbracefill\crcr\noalign{\kern3\p@\nointerlineskip}
      $\hfil\displaystyle{#1}\hfil$\crcr}}}\limits}
 \def\underbrace#1{\mathop{\vtop{\m@th\ialign{##\crcr
      $\hfil\displaystyle{#1}\hfil$\crcr\noalign{\kern3\p@\nointerlineskip}
      \upbracefill\crcr\noalign{\kern3\p@}}}}\limits}

  \let\@ldmatrix\matrix
  \let\end@ldmatrix\endmatrix
  \def\matrix{\sam\@ldmatrix}
  \def\endmatrix{\end@ldmatrix\mas}
  \let\@ldgather\gather
  \let\end@ldgather\endgather
  \def\gather{\sam\@ldgather}
  \def\endgather{\end@ldgather\mas}
  \let\@ldalign\align
  \let\end@ldalign\endalign
  \def\align{\sam\@ldalign}
  \def\endalign{\end@ldalign\mas}
  \let\@ldaligned\aligned
  \let\end@ldaligned\endaligned
  \def\aligned{\sam\@ldaligned}
  \def\endaligned{\end@ldaligned\mas}
  \let\@ldtag\tag
  \def\tag{\sam\@ldtag}
   %

   \let\MinCDArrowWidth\minCDaw@




\newskip\insertskipamount\newskip\inserthardskipamount
\insertskipamount 6pt plus2pt 
\inserthardskipamount 6pt
\def\insertskip{\vskip\insertskipamount}
\newcount\SplitTest
\def\SetSplitTest{\SplitTest\insertpenalties
  \insert\topins{\floatingpenalty1}%
  \advance\SplitTest-\insertpenalties}
\def\midinsert{\par
 \SaveLastSkip\penalty-150\SetSplitTest\RestoreLastSkip
 \ifnum\SplitTest=-1
  \@midfalse\p@gefalse\else\@midtrue\fi\@ins}
\def\@ins{\par\begingroup\setbox\z@\vbox\bgroup%
  \vglue\inserthardskipamount}
\def\endinsert{\egroup 
  \if@mid \dimen@\ht\z@ \advance\dimen@\dp\z@
    \advance\dimen@\insertskipamount
    \advance\dimen@\pagetotal\advance\dimen@-\pageshrink
    \ifdim\dimen@>\pagegoal\@midfalse\p@gefalse\fi\fi
  \if@mid%
    \ifdim\lastskip<\insertskipamount\removelastskip\insertskip\fi
    \nointerlineskip\box\z@\penalty-200\insertskip
  \else%
    \SaveLastSkip
    \insert\topins{\penalty100 
    \splittopskip\z@skip
    \splitmaxdepth\maxdimen \floatingpenalty\z@
    \ifp@ge \dimen@\dp\z@
    \vbox to\vsize{\unvbox\z@\kern-\dimen@}
    \else \box\z@\nobreak\insertskip\fi}
    \RestoreLastSkip
   \fi\endgroup}


  \newcount\notenumber
  
  \def\note{\advance\notenumber by 1
    \footnote{\the\notenumber)}}

  \newbox\footbox

  \def\footnote#1{\let\@sf\empty
    \ifhmode\edef\@sf{\spacefactor\the\spacefactor}\/\fi
    \sam${}^{\fam0 #1}$\@sf\vfootnote{#1}}%

  \def\vfootnote#1{\insert\footins\bgroup
     \interlinepenalty100 \splittopskip=1pt
     \floatingpenalty=20000
     \leftskip=0pt\rightskip=0pt%
     \parindent=.3em
     \Smallfonts\rm
     \FootItem@{#1}
     \futurelet\next\fo@t}

  \def\FootItem@#1{\par\hangafter1\hangindent=\FootHang
     \setbox0=\hbox{\ignorespaces#1\unskip}%
     \dimen0=.4em\SetOverhang@
     \noindent\rlap{\box0}\kern\Overhang\ignorespaces}


  \def\fo@t{\ifcat\bgroup\noexpand\next \let\next\f@@t
    \else\let\next\f@t\fi \next}
  \def\f@@t{\bgroup\aftergroup\@foot\let\next}
  \def\f@t#1{\baselineskip=10pt\lineskip=1pt
            \lineskiplimit=0pt #1\@foot}%
  \def\@foot{
        \hbox{\vrule height0pt depth5pt width0pt}
        \egroup}
  \skip\footins=12 pt plus 0pt minus 0pt 
  \count\footins=1000 
  \dimen\footins=8in 



 \def\osumess#1{\EdSpider{\immediate\write16{Line \the\inputlineno: #1}}}%
 \def\HideEdStuff{\gdef\EdSpider##1{}}

 \font\BigSym=cmmi10 scaled \magstep 4

 \def\change{\InLMargin{\hbox{\BigSym \char63\kern10pt}}}

 \def\beginchange{\InLMargin{\hbox{\sam\twelvepoint$\heartsuit$\kern10pt}}}

 \def\endchange{\InLMargin{\hbox{\sam\twelvepoint$\spadesuit$\kern10pt}}}

 \def\InLMargin#1{\strut\vadjust{%
     \kern-\strutdepth
     \vtop to \strutdepth{%
         \baselineskip\strutdepth
         \llap{\sam$\smash{\hbox{\EdSpider{#1}}}$}\null}}}

 \def\strutdepth{\dp\strutbox}
 \def\strutheight{\ht\strutbox}

 \def\NoteInRMargin#1{\strut\vadjust{%
     \kern-1.001\strutdepth
     \vtop to \strutdepth{%
       \baselineskip\strutdepth
       \vss\rlap{\ninepoint\unskip\hskip\hsize
         \vtop to 0pt{%
           \hsize=16em\hfuzz=\hsize
           \leftskip=10pt%
           \rightskip=0pt plus 10000pt%
           \baselineskip=9.8pt\lineskip=.2pt%
           \let\\\break
           \noindent\EdSpider{#1}\vss}%
                \kern10pt}\hbox{}}
       }}

 \def\ednote#1{\NoteInRMargin{\tentt #1}}

 \def\cbar{\InLMargin{%
      \dimen0=\strutdepth\advance\dimen0 by \lineskip
      \vrule width 3pt
      height \strutheight depth \dimen0 \kern
      3pt}}

 \def\ccbar{\InLMargin{%
      \dimen0=2\strutdepth\advance\dimen0 by 2\lineskip
      \vrule width 3pt
        height 3\strutheight depth \dimen0 \kern
      3pt}}

 \newinsert\TRMargIns
 \dimen\TRMargIns=\maxdimen

  \def\Ednote#1{\insert\TRMargIns{%
       \vbox to 0pt{\hsize=140pt\hfuzz=\hsize
           \leftskip=6pt%
           \rightskip=0pt plus 10000pt%
           \baselineskip=9.8pt\lineskip=.2pt%
           \let\\\break
           \SetPageRemainder
           \vglue540pt\vglue-\PageRemainder
           \noindent\EdSpider{\tentt #1}\vss}%
       \smallskip}}

 \def\KillEdStuff{\def\ednote##1{}\def\Ednote##1{}%
      \let\change\relax\let\beginchange\relax\let\endchange\relax
       \let\cbar\relax\let\ccbar\relax}


  \topskip=12pt
  \newskip\StdBaselineskip 
  \StdBaselineskip 12pt
  \lineskip=1.1pt
  \lineskiplimit=.8pt
  \widowpenalty=10000 
  \clubpenalty=10000  
  \abovedisplayskip=6pt plus 1pt minus 1pt
  \abovedisplayshortskip=3pt plus 1.5pt
  \belowdisplayskip=6pt plus 1pt minus 1pt
  \belowdisplayshortskip=5pt plus 1pt minus 1pt
  \hfuzz=1.5pt   

  \def\StdPretolerance{100}
  \tolerance=\StdPretolerance

  \newdimen\StdMathsurround
  \StdMathsurround=1.5pt 
  \mathsurround=\StdMathsurround
  \Mas                   

   \def\prose{\relax\hbox{\kern.6\StdMathsurround}}
  
  \def\StdParskip{0pt}    
  \parskip=\StdParskip
  \parindent=0.5cm
 

  \def\Times{ptmr  } 
  \def\TimesI{ptmri  } 
  \def\TimesB{ptmb  }
  \def\TimesBI{ptmbi  }
  \def\HelveticaN{phvrrn }

  =\Times at 10bp
  =\TimesB at 10bp
  \font\tenit=\TimesI at 10bp
  =\TimesBI at 10bp

  \font\tenmrm=cmr10  


    =\Times at 9bp 
    \font\nineit=\TimesI at 9bp 
    =\TimesB at 9bp 
    =\TimesBI at 9bp 

    =\HelveticaN at 9bp 


  =\Times at 12bp
  \font\twelveit=\TimesI at 12bp
  =\TimesB at 12bp


  \font\titleit=\TimesI at 14.4bp
  =\TimesB at 14.4bp

 \SetAuthorHead{AuthorHead} 
 \SetTitleHead{TitleHead}  


  \def\lBr{\raise.125ex\hbox{[\kern.1125ex}}
  \def\rBr{\raise.125ex\hbox{\kern.1125ex]}}

 \setbox\footbox=\hbox{\Smallfonts 2)~}



  \bgroup
  \catcode`\@=11 
  \gdef\itSpacing{%
     \xspaceskip=.31em plus.1em minus.05em \sfcode `f=2001
     \itWarning@\let\itWarning@\itWarning@@}
  \gdef\itSpacingOff{%
     \xspaceskip=0pt \sfcode `f=1000
     \let\itWarning@\relax}
   \global\let\itWarning@\relax
  \gdef\itWarning@@{\errmessage{%
  Special italic spacing already in force
  (you have probably omitted an ``endth'').
  See itSpacing macro in osuPSfnt.sty
         }}
  \egroup

 \fontdimen1\titlebf=0.0pt
 \fontdimen2\titlebf=3.6135pt
 \fontdimen3\titlebf=2.8908pt
 \fontdimen4\titlebf=1.44539pt
 \fontdimen5\titlebf=6.64882pt
 \fontdimen6\titlebf=14.45398pt
 \fontdimen7\titlebf=1.60439pt

 \fontdimen1\tenbi=0.26794pt
 \fontdimen2\tenbi=2.50937pt
 \fontdimen3\tenbi=2.00749pt
 \fontdimen4\tenbi=1.00374pt
 \fontdimen5\tenbi=4.59717pt
 \fontdimen6\tenbi=10.03749pt
 \fontdimen7\tenbi=1.11415pt

 \fontdimen1\twelverm=0.0pt
 \fontdimen2\twelverm=3.01125pt
 \fontdimen3\twelverm=2.409pt
 \fontdimen4\twelverm=1.2045pt
 \fontdimen5\twelverm=5.39615pt
 \fontdimen6\twelverm=12.045pt
 \fontdimen7\twelverm=1.33699pt

 \fontdimen1\twelveit=0.27731pt
 \fontdimen2\twelveit=3.01125pt
 \fontdimen3\twelveit=2.409pt
 \fontdimen4\twelveit=1.2045pt
 \fontdimen5\twelveit=5.37207pt
 \fontdimen6\twelveit=12.045pt
 \fontdimen7\twelveit=1.33699pt

 \fontdimen1\twelvebf=0.0pt
 \fontdimen2\twelvebf=3.01125pt
 \fontdimen3\twelvebf=2.409pt
 \fontdimen4\twelvebf=1.2045pt
 \fontdimen5\twelvebf=5.5407pt
 \fontdimen6\twelvebf=12.045pt
 \fontdimen7\twelvebf=1.33699pt

 \fontdimen1\tenrm=0.0pt
 \fontdimen2\tenrm=2.50937pt
 \fontdimen3\tenrm=2.00749pt
 \fontdimen4\tenrm=1.00374pt
 \fontdimen5\tenrm=4.49678pt
 \fontdimen6\tenrm=10.03749pt
 \fontdimen7\tenrm=1.11415pt

 \fontdimen1\tenit=0.27731pt
 \fontdimen2\tenit=2.50937pt
 \fontdimen3\tenit=2.00749pt
 \fontdimen4\tenit=1.00374pt
 \fontdimen5\tenit=4.47672pt
 \fontdimen6\tenit=10.03749pt
 \fontdimen7\tenit=1.11415pt

 \fontdimen1\tenbf=0.0pt
 \fontdimen2\tenbf=2.50937pt
 \fontdimen3\tenbf=2.00749pt
 \fontdimen4\tenbf=1.00374pt
 \fontdimen5\tenbf=4.61723pt
 \fontdimen6\tenbf=10.03749pt
 \fontdimen7\tenbf=1.11415pt

 \fontdimen1\ninerm=0.0pt
 \fontdimen2\ninerm=2.25842pt
 \fontdimen3\ninerm=1.80673pt
 \fontdimen4\ninerm=0.90337pt
 \fontdimen5\ninerm=4.0471pt
 \fontdimen6\ninerm=9.03374pt
 \fontdimen7\ninerm=1.00273pt

 \fontdimen1\nineit=0.27731pt
 \fontdimen2\nineit=2.25842pt
 \fontdimen3\nineit=1.80673pt
 \fontdimen4\nineit=0.90337pt
 \fontdimen5\nineit=4.02904pt
 \fontdimen6\nineit=9.03374pt
 \fontdimen7\nineit=1.00273pt

 \fontdimen1\ninebf=0.0pt
 \fontdimen2\ninebf=2.25842pt
 \fontdimen3\ninebf=1.80673pt
 \fontdimen4\ninebf=0.90337pt
 \fontdimen5\ninebf=4.15552pt
 \fontdimen6\ninebf=9.03374pt
 \fontdimen7\ninebf=1.00273pt


 \newcount\MaxSpaceFactor
 \MaxSpaceFactor=3000 

 \def\ItemStyle{\rm}
 \def\NrStyle{\rm}
 \def\ItemItemStyle{\rm}

 \MaxItemTag{(iii)}
 \MaxItemItemTag{(iii)}
 \MaxNrTag{(2)}
 \MaxFootTag{2)}
 \def\ReferenceHang{30pt}

 \catcode`\@=\active


\loadbold

=\Times  
=\Times scaled750
=\Times scaled650
\font\rms=\Times scaled 920 

=\TimesBI scaled 860
=\TimesI scaled 860

\textfont0=\rrm  
\scriptfont0=\erm 
\scriptscriptfont0=\srm

\def\Augment#1#2{%
    \toks0\expandafter{#1}\toks2{#2}%
    \edef#1{\the\toks0\the\toks2}}

 \font\twelverma=\Times  scaled 1200
 \font\tenrma=\Times  scaled 1000
 \font\ninerma=\Times scaled 920
 =\Times scaled 840
 \font\sevenrma=\Times scaled 760
 =\Times scaled 680
 \font\fiverma=\Times scaled 600

 \Augment\tenpoint{%
  \textfont0=\tenrma  \scriptfont0=\sevenrma  
  \scriptscriptfont0=\fiverma  }

 \Augment\ninepoint{%
  \textfont0=\ninerma  \scriptfont0=\sevenrma 
  \scriptscriptfont0=\fiverma}

 \Augment\twelvepoint{%
  \textfont0=\twelverma  \scriptfont0=\ninerma  
  \scriptscriptfont0=\sevenrma}

\mathsurround=1pt
\hsize=13.45truecm
\vsize=19.5truecm
\hoffset=1.25truecm
\voffset=2truecm
\advance\baselineskip by 2pt

\predefine\til{\~}
\def\~#1{\relax\ifmmode\widetilde{#1}\else\til{#1}\fi}

\redefine \le{\leqslant}
\redefine \ge{\geqslant}
\define \wt#1{\mathaccent"0365{#1}}
\define \wh#1{\mathaccent"0362{#1}}

\def\CO{{O}}
\def\CM{{M}}

\define\cdvf{complete discrete valuation field }

\define \iss{\,\Mathaccent{\raise -.8 ex\hbox{$\widetilde{}$\kern.1em}}\rightarrow\,}

\define \minn{\operatorname{\fam0 min\,}}

\define \innf{\operatorname{\fam0 inf\,}}

\define \chr{\mathop{\fam0 char}\,}

\Mas
\HideEdStuff
\rm 
 

\def\issn{{\nineit ISSN 1464-8997 (on line) 1464-8989 (printed)}}

\def\gtp{{\nineit Published 10 December 2000: \ \copyright\ Geometry \& 
Topology Publications}}

\def\gtv3{{\nineit Geometry \& Topology Monographs, Volume 3 (2000) --
Invitation to higher local fields}}


\def\lione
{{\rms Geometry \& Topology Monographs}}

\def \litwo{{\rms Volume 3: Invitation to higher local fields
}} 

\def\tinfo #1.#2.#3-#4
{{
\noindent  {\lione} \hfill 
\par 
\vskip-1.5pt
\noindent {\litwo} \hfill
\par 
\vskip-1,5pt
\noindent {\rms Part #1, section #2, pages #3--#4} \hfill
\vskip24pt 
}}

\def\tinfos #1.#2.#3-#4
{{
\noindent  {\lione} \hfill 
\par 
\vskip-1.5pt
\noindent {\litwo} \hfill
\par 
\vskip-1.5pt
\noindent {\rms Pages #3--#4} \hfill
\vskip24pt 
}}

\def\tinfoi #1
{{
\noindent  {\lione} \hfill 
\par 
\vskip-1.5pt
\noindent {\litwo} \hfill
\par 
\vskip-1.5pt
\noindent {\rms Pages iii--xi: Introduction and contents} \hfill
\vskip26pt 
}}


  \def\titlepagehead{\hfil}

  \newif\iftitlepage\titlepagefalse
  \newif\ifblankpage\blankpagefalse
  \def\makeheadline{
     \ifblankpage{}\else%
     \iftitlepage
\vbox{\line{\vbox to 8.5pt{}
\ninerm
\copy\HLinebox \hfill
\hglue5mm\ninebf\folio 
\titlepagehead}}%
      \else
\vbox{\ifodd\pageno\rightheadline\else\leftheadline\fi}%
      \fi\vskip 12pt\fi}%
     \def\rightheadline{\line{\vbox to 8.5pt{}%
      \ninerm
\copy\TitleBox \hfill
\hglue5mm\ninebf\folio}}%
     \def\leftheadline{\line{\vbox to 8.5pt{}%
        \unskip\ninerm\unskip\ninebf\folio\hglue5mm
 \hfill \copy\AuthorBox
}}

 \footline={\ifblankpage{}\else
\iftitlepage\ninepoint\sam\hfill
\line{\vbox to 8.5pt{}
\copy\TFLinebox
\hfill
\hglue5mm 
}
            \else
\ninepoint\sam\hfill
\line{\vbox to 8.5pt{}
\copy\FLinebox
\hfill 
\hglue5mm
}
\hfil\fi\global\titlepagefalse\fi}

\def\blankpage{{\blankpagetrue\noindent\hbox to 10pt{\hss}\vfill
\pagebreak}}

\tenpoint\rm 
 